\documentclass[12pt]{article}
\usepackage{mathtools}
\usepackage{appendix}

\usepackage{amsmath,latexsym,amsfonts,enumerate}
\usepackage{amssymb,amsthm}
\usepackage{calrsfs} 
\DeclareMathAlphabet{\pazocal}{OMS}{zplm}{m}{n}
\usepackage{tikz-cd}
\usepackage{tikz}
\usepackage{comment}

\usepackage[LGR,T1]{fontenc} 
\usepackage[utf8]{inputenc} 
\usepackage{lmodern} 
\usepackage[greek,english]{babel} 
\usepackage{alphabeta} 
\usepackage{indentfirst} 
\usepackage{biblatex}
\usepackage{csquotes}
\usepackage[colorlinks=true,urlcolor=black,linkcolor=blue,citecolor=blue]{hyperref}
\usepackage{cleveref}
\usepackage{pgfplots}
\usepackage{tikz-3dplot}
\usepackage{faktor}
\usepackage[colorlinks=true,urlcolor=black,linkcolor=blue,citecolor=blue]{hyperref}

\usepackage[T3,T1]{fontenc}
\DeclareSymbolFont{tipa}{T3}{cmr}{m}{n}
\DeclareMathAccent{\invbreve}{\mathalpha}{tipa}{16}

\tdplotsetmaincoords{60}{115}
\pgfplotsset{compat=newest}

\newtheorem{theorem}{Theorem}[section]
\newtheorem*{theorem*}{Theorem}
\newtheorem{lemma}{Lemma}[section]
\newtheorem*{lemma*}{Lemma}

\theoremstyle{definition}
\newtheorem{definition}{Definition}[section]
\newtheorem*{remark*}{Remark}
\newtheorem{remark}{Remark}

\newcommand{\ignore}[1]{}

\title{Classification of analytic $\text{SO}^o(p,q)$-actions on closed $(p+q-1)$-dimensional manifolds II :  $p \geq 3$, $q=1,2$}
\author{{\bf Spyridon Lentas} }
\date{}

\begin{document}

\maketitle

\begin{abstract}
	This paper provides a classification of analytic actions of the semi-orthogonal groups $\text{SO}^\circ(p,1)$ and $\text{SO}^\circ(p,2)$, for $p \geq 3$, on closed, connected manifolds of dimension $p$ and $p+1$ respectively. It is a sequel to a previous paper, which treats the case $p,q \geq 3$. Adapting Uchida's methods, we explicitly construct analytic actions of $\text{SO}^\circ(p,1)$ on $\text{S}^{p}$ and $\text{S}^{p-1} \times \text{S}^{1}$, and analytic actions of $\text{SO}^\circ(p,2)$ on $\text{S}^{p+1}$, $\text{S}^{p} \times \text{S}^{1}$ and $\text{S}^{p-1} \times \text{S}^{2}$, and we classify them up to analytic isomorphism, using a result of Hitchin and Schneider's classification of analytic $\text{SL}_2(\mathbb{R})$ actions on surfaces. The main result demonstrates that any analytic $\text{SO}^\circ(p,q)$ action, $q=1,2$, on a closed, connected $(p+q-1)$-dimensional manifold is covered by one of the constructed actions or by an action on $\text{SO}^\circ(p,q) \times_P \text{S}^1$, where $P$ is a maximal parabolic subgroup of $\text{SO}^\circ(p,q)$.

\end{abstract}

\section{Introduction} \label{intro}
This is the second of two papers, both based on the author's PhD dissertation \cite{Lentas_thesis}, on analytic actions of the semi-orthogonal groups $\text{SO}^\circ(p,q)$, $p \geq 3$, on closed, connected manifolds of dimension $p+q-1$. The first paper \cite{Lentas_SOpq_I} deals with the case $p,q \geq 3$; here we complete the picture by treating the groups $\text{SO}^\circ(p,1)$ and $\text{SO}^\circ(p,2)$. In other words, we study analytic $\text{SO}^\circ(p,1)$ actions on manifolds of dimension $p$ and analytic $\text{SO}^\circ(p,2)$ actions on manifolds of dimension $p+1$.

The dimension $p+q-1$ is the first one in which interesting actions appear. Indeed, $\text{SO}^\circ(p,q)$ admits no nontrivial action on a closed, connected manifold of dimension smaller than $p+q-2$, while in dimension $p+q-2$ the only action, essentially, is the one on the homogeneous space $\text{SO}^\circ(p,q) \big/ P_{\text{null}}$. Here, and throughout the paper, $P_{\text{null}} \leq \text{SO}^\circ(p,q)$ denotes the maximal parabolic subgroup that stabilises an isotropic line in $\mathbb{R}^{p+q}$, where $\mathbb{R}^{p+q}$ is the standard representation of $\text{SO}^\circ(p,q)$ and carries a scalar product of signature $(p,q)$. Similarly, $G_{\text{null}} \leq \text{SO}^\circ(p,q)$ denotes the stabiliser of a null vector, and, following \cite{Lentas_SOpq_I}, an orbit diffeomorphic to $\text{SO}^\circ(p,q) \big/ G_{\text{null}}$ is called a \textit{nullcone orbit}; as a topological space, it is a component of the nullcone of the scalar form of signature $(p,q)$ on $\mathbb{R}^{p+q}$.

Our approach has two parts. First, we single out a small family of model spaces, namely $\text{S}^p$ and $\text{S}^{p-1} \times \text{S}^1$ when $q=1$, and $\text{S}^{p+1}$, $\text{S}^p \times \text{S}^1$ and $\text{S}^{p-1} \times \text{S}^2$ when $q=2$, and we classify, up to analytic isomorphism, the analytic actions on them for which a maximal compact subgroup acts in a standard way. Second, we prove that every closed, connected $(p+q-1)$-dimensional manifold with a nontrivial analytic $\text{SO}^\circ(p,q)$ action is $\text{SO}^\circ(p,q)$-equivariantly covered either by one of these model spaces or by the associated bundle $\text{SO}^\circ(p,q) \times_{P_{\text{null}}} \text{S}^1$, on which $P_{\text{null}}$ acts on the circle through an analytic flow, see \cite{Lentas_SOpq_I}. At the level of diffeomorphism types, this yields the following two theorems.

\begin{theorem} \label{general thm (p,1)}
Let $p \geq 3$ and let $M$ be a closed, connected $p$-dimensional manifold with an analytic action of $\textup{SO}^\circ(p,1)$.
\begin{itemize}
    \item If some point of $M$ is fixed by $\textup{SO}(p)$, then $\textup{S}^p$ covers $M$.
    \item If no point of $M$ is fixed by $\textup{SO}(p)$ and $M$ contains no nullcone orbit, then $\textup{S}^{p-1} \times \textup{S}^1$ covers $M$.
    \item If no point of $M$ is fixed by $\textup{SO}(p)$ and $M$ contains a nullcone orbit, then $\textup{SO}^\circ(p,1) \times_{P_{\textup{null}}} \textup{S}^1$ covers $M$, where, writing $P_{\textup{null}}=M_P A_P N_P$ for the Langlands decomposition, $P_{\textup{null}}$ acts on $\textup{S}^1$ by a flow via $A_P$.
\end{itemize}
\end{theorem}

\begin{theorem} \label{general thm (p,2)}
Let $p \geq 3$ and let $M$ be a closed, connected $(p+1)$-dimensional manifold with an analytic action of $\textup{SO}^\circ(p,2)$. Denote by $\textup{SO}(p)$ the subgroup $\textup{SO}(p) \times \{1\}$ of the maximal compact subgroup $\textup{SO}(p) \times \textup{SO}(2) \leq \textup{SO}^\circ(p,2)$, and define $\textup{SO}(2)$ analogously.
\begin{itemize}
    \item If $\textup{SO}(p)$ and $\textup{SO}(2)$ both have fixed points, then $\textup{S}^{p+1}$ covers $M$.
    \item If $\textup{SO}(p)$ has fixed points but $\textup{SO}(2)$ does not, then $\textup{S}^p \times \textup{S}^{1}$ covers $M$.
    \item If $\textup{SO}(2)$ has fixed points but $\textup{SO}(p)$ does not, then $\textup{S}^{p-1} \times \textup{S}^{2}$ covers $M$.
    \item If neither $\textup{SO}(p)$ nor $\textup{SO}(2)$ has fixed points, then $\textup{SO}^\circ(p,2) \times_{P_{\textup{null}}} \textup{S}^1$ covers $M$, where, writing $P_{\textup{null}}=M_P A_P N_P$ for the Langlands decomposition, $P_{\textup{null}}$ acts on $\textup{S}^1$ by a flow via $A_P$.
\end{itemize}
When $p=3$, the only exception is a transitive action covered by $\textup{SO}^\circ(3,2) \big/ P_{\textup{min}}$, with $P_{\textup{min}} \leq \textup{SO}^\circ(3,2)$ a minimal parabolic subgroup.
\end{theorem}

 The classification theorems in Sections \ref{SO(p,1) axns} and \ref{SO(p,2) axns} are more precise: they also identify the action on the covering space as one coming from the basic constructions introduced in those sections. Although the overall structure parallels the case $p,q \geq 3$ of \cite{Lentas_SOpq_I}, two new phenomena appear. For $q=1$, $M$ may contain more than two disjoint orbits diffeomorphic to $\faktor{\text{SO}^\circ(p,1)}{ \text{SO}^\circ(p-1,1)}$, and knowing that $\text{SO}(p)$ has no fixed points does not by itself determine the topology of $M$. In addition, a basic ingredient for the results in the cases where $q\neq 2$ is the fixed point set of a specific subgroup, $H \simeq \text{SO}(p-1)$, of $G$, see the Notation section. For $q=2$, the fixed point set of $H$ is a surface rather than a union of circles, which was the case for $q\neq2$, so the 1-dimensional methods of \cite{Lentas_SOpq_I} have to be replaced by 2-dimensional ones.

 The paper is organised as follows. Section \ref{SO(p,1) axns} is devoted to $q=1$. In this case the fixed point set of $H$ is still a finite union of circles. 
 We construct analytic $\text{SO}^\circ(p,1)$ actions on $\text{S}^p$, see also \cite{Uchida_SOpqII}, and on $\text{S}^{p-1} \times \text{S}^1$, extending the standard action of the maximal compact subgroup $\text{SO}(p)$. The actions on $\text{S}^{p-1} \times \text{S}^1$ are built from what we call \textit{simple flows} on $\text{S}^1$; by a result of \cite{HITCHIN1991359}, these flows are determined up to analytic isomorphism by their number of fixed points, the Jacobians at the fixed points, and a \textit{global invariant} $\mu$. The section ends by showing that any analytic $\text{SO}^\circ(p,1)$ action on a closed, connected $p$-dimensional manifold is equivariantly covered by one of these actions or, when a nullcone orbit is present, by $\text{SO}^\circ(p,1) \times_{P_{\text{null}}} \text{S}^1$.

 Section \ref{SO(p,2) axns} treats $q=2$. Building on \cite{Uchida_SOpqII}, we associate to every analytic $\text{SO}^\circ(p,2)$ action on $\text{S}^p \times \text{S}^1$, $\text{S}^{p-1} \times \text{S}^2$ or $\text{S}^{p+1}$ that extends the standard $\text{SO}(p) \times \text{SO}(2)$ action an analytic action on the 2-dimensional fixed point set of $H \simeq \text{SO}(p-1)$: an $\text{SO}^\circ(1,2)^+$ (see (\ref{SO(1,2)+ def})) action on $\mathbb{T}^2$, an $\text{SO}^\circ(1,2)$ action on $\text{S}^2$, and an $\text{SO}^\circ(1,2)^+$ action on $\text{S}^2$, respectively. In Section \ref{(p,2) case constr} we reverse this process and build the $\text{SO}^\circ(p,2)$ actions out of these data, see Sections \ref{(p,2)-basic const I}, \ref{(p,2)-basic const II} and \ref{(p,2)-basic const III}. The isomorphism $\text{SO}^\circ(1,2) \simeq \text{PSL}_2(\mathbb{R})$ then allows us to classify the resulting actions using Schneider's classification of analytic $\text{SL}_2(\mathbb{R})$ actions on surfaces \cite{SchneiderSL2axns}. We conclude the section by proving that every analytic $\text{SO}^\circ(p,2)$ action on a closed, connected $(p+1)$-dimensional manifold is equivariantly covered by one of the constructed actions.
 
Many results in this paper have proofs which are completely analogous to proofs of results in \cite{Lentas_SOpq_I}. In order to keep the text concise and readable, we omit those proofs and refer to the corresponding ones in \cite{Lentas_SOpq_I}.

\subsection{Background and Motivation}

 The starting point of this work is Uchida's study \cite{Uchida} of smooth $\text{SO}^\circ(p,q)$ actions on $\text{S}^{p+q-1}$ for which a maximal compact subgroup acts on the sphere by the standard orthogonal action. Uchida encoded each such action by a pair $(\Phi, f)$, consisting of a smooth flow $\Phi$ on $\text{S}^1$ and a smooth function $f : \text{S}^1 \rightarrow \mathbb{R}\text{P}^1$ subject to certain compatibility conditions, and showed that this gives a one to one correspondence. The circle in question is the fixed point set of $H$. Note that $H$ is a principal isotropy group whose orbits under the maximal compact subgroup have codimension 1. The value of $f$ at a point records the isotropy algebra of the $\text{SO}^\circ(p,q)$ action there, see also \cite{Lentas_SOpq_I}. However, Uchida's correspondence classifies the actions only up to homeomorphism.

 Such pairs of a flow and a function first appeared in Asoh's work, see \cite{AsohSL2C}, on $\text{SL}_2(\mathbb{C})$ actions on $\text{S}^3$, were adapted by Uchida in \cite{Uchida} to the setting described above, and were further refined by Mukoyama in \cite{MukoyamaSp2r} for smooth $\text{Sp}(2,\mathbb{R})$ actions on $\text{S}^4$. Mukoyama's version was in turn used by Uchida for smooth $\text{SO}^\circ(p,2)$ actions on $\text{S}^{p+1}$ in \cite{Uchida_SOpqII} and for $\text{Sp}(p,q)$ actions on $\text{S}^{4p+4q-1}$, see \cite{UchidaSppq}, and by Mukoyama for smooth $\text{SU}(p,q)$ actions on $\text{S}^{2p+2q-1}$ and on $\mathbb{R} \text{P}^{p+q-1}$ in \cite{MukoyamaSUpq}. Recall that $\text{SL}_2(\mathbb{R})$, $\text{SL}_2(\mathbb{C})$ and $\text{Sp}(2,\mathbb{R})$ are locally isomorphic to $\text{SO}^\circ(1,2)$, $\text{SO}^\circ(1,3)$ and $\text{SO}^\circ(2,3)$ respectively.

 In the analytic setting, Schneider classified the $\text{SL}_2(\mathbb{R})$ actions on compact surfaces in  \cite{SchneiderSL2axns}, and this classification plays a central role in Section \ref{SO(p,2) axns}. His actions fall into three families: actions all of whose orbits are closed circles; effective actions, which are classified by a continuous invariant; and non-effective actions, which are classified by finitely many discrete invariants and may have several open orbits separated by closed circle orbits. The last family resembles both Asoh's analytic actions and the analytic $\text{SO}^\circ(p,1)$ actions of Section \ref{SO(p,1) axns}. By contrast, an analytic $\text{SO}^\circ(p,q)$ action with $p \geq 3$ and $q \geq 2$ has either at most two open orbits, or arbitrarily many open orbits, all of them nullcone orbits. In the smooth setting of Uchida's result on $\text{S}^{p+q-1}$, there may be more than two open orbits, all but two of which are nullcone orbits. See also \cite{Stowe_SL2R_surfaces} for a treatment of the analytic $\text{SL}_2(\mathbb{R})$ actions on surfaces and \cite{Gorbatsevich_loc_trans_surfaces} for a treatment of analytic Lie group actions on surfaces.

 Related classification results exist for $\text{SL}_n(\mathbb{R})$. Uchida classified the analytic $\text{SL}_n(\mathbb{R})$ actions on $\text{S}^n$ in \cite{Uchida_clssf_rlanalytic_SLn}, and Fisher and Melnick later classified both the analytic and the smooth $\text{SL}_n(\mathbb{R})$ actions on closed $n$-dimensional manifolds in \cite{FisherMelnick}. Son extended the analytic classification to manifolds of dimension $m$ with $n\leq m \leq 2n-3$, and classified the smooth fixed point free actions in the same range of dimensions, see \cite{Son_Reanalytic_SLn_n<m<2n-2}.

\subsection{Notation}
We fix the following notation for the rest of the paper.

\begin{itemize}
    \item[-] For a Lie group $A$, its identity component is denoted by $A^\circ$.
    \item[-] $ \displaystyle G := \text{SO}^\circ(p,q)= \left\{ X \in \text{SL}_{p+q}(\mathbb{R}) : X \, I_{p,q} \, X^T = I_{p,q} \right\}$, $p\geq3$,\\ where $I_{p,q} = \begin{bmatrix}
    -I_p & \\
    & I_q
\end{bmatrix}$
and $I_p$, $I_q$ denote the $p \times p$ and $q \times q$ identity matrices.
    \item[-] $ \displaystyle K := \left\{ \begin{bmatrix}
    A & \\ & B
\end{bmatrix} \in G : A \in \text{SO}(p) , B \in \text{SO}(q) \right\} \leq G   $, so that $K \simeq \text{SO}(p) \times \text{SO}(q)$.
    \item[-] $ \displaystyle H := \left\{ \begin{bmatrix}
    1 & && \\
    & \Tilde{A} && \\
    & && I_q  
\end{bmatrix} \in K : \Tilde{A} \in \text{SO}(p-1) \right\} \leq K$,\\
so that $H \simeq \text{SO}(p-1)$.
    \item[-] $\mathcal{M}(p,q) = \left\{ \begin{bmatrix}
    \text{cosh}(\theta) & & \text{sinh}(\theta) & \\
    & I_{p-1} && \\
    \text{sinh}(\theta) & & \text{cosh}(\theta) & \\
    &&& I_{q-1}
\end{bmatrix} \right\} \leq G$
    \item[-] $j_1 = \begin{bmatrix}
        -1&& \\ & -1 & \\ && I_{p+q-2}
    \end{bmatrix} \in \text{SO}(p) \times \{ I_q \} \leq K$
    \item[-] $j_2 = \begin{bmatrix}
        I_p &&& \\ & -1 && \\ &&-1 & \\ &&& I_{q-1}
    \end{bmatrix} \in \{ I_p \} \times \text{SO}(q) \leq K$, defined only when $q\geq2$.
    \item[-] When a group $G$ acts on a space $X$, we write $g \star x$ for the image of $x \in X$ under $g \in G$, except when a matrix $A$ acts linearly on a vector $v$, in which case we simply write $Av$.
    \item[-] For a Lie group $A$, a Lie subgroup $B \leq A$ and a $B$-space $X$, we denote by $\displaystyle A \times_B X$ the quotient of $A \times X$ by the equivalence relation $\displaystyle (a,x) \sim (ab, b^{-1}x)$. The group $A$ acts on $A \times_B X$ by left multiplication, $\displaystyle \Tilde{a} \star (a,x) = (\Tilde{a}a,x)$.

\end{itemize}

\subsection*{Acknowledgments}
The author is deeply grateful to his advisor, Prof. Karin Melnick, for her guidance and constant support during this project, and would like to thank Prof. Bill Goldman for many helpful conversations.

\section{Classification of analytic $\text{SO}^\circ(p,1)$, $p \geq 3$, actions} \label{SO(p,1) axns}

Here we consider the case $q=1$, $p \geq 3$, namely an analytic action of $G= \text{SO}^\circ(p,1)$ on a closed, connected manifold $M$ of dimension $p$. We see that these actions are covered by actions of $G$ on either $\text{S}^p$ or $\text{S}^{p-1} \times \text{S}^1$. While the topology of these spaces is analogous to the corresponding spaces of the $p,q \geq 3$ case, see \cite{Lentas_SOpq_I}, the geometry can differ. The fixed point set of $H = \text{SO}(p-1)$ is again going to be a union of circles. In $\text{SO}(p)$ there exists the involution $j_1$ which forces the fixed points of $\text{SO}(p)$ on one of these circles to be at most 2. However, in contrast to the $q\geq 3$ case, there is no longer a $j_2 = \text{diag}(I_p,-I_2,I_{q-2})$. That means there could be more than two points whose isotropy group is the standard copy of $\text{SO}^\circ(p-1,1)$, and so more than two disjoint orbits diffeomorphic to
\[  \faktor{\text{SO}^\circ(p,1)}{ \text{SO}^\circ(p-1,1)}  \]
Another notable difference is that absence of $K$ fixed points does not imply the existence of nullcone orbits. The actions of $K = \text{SO}(p)$ on $\text{S}^p$ and $\text{S}^{p-1} \times \text{S}^1$ are via rotations with respect to an axis in the former case and the standard orthogonal action on the first factor in the latter.

Furthermore, a result analogous to \cite[Lemma 2.2]{Lentas_SOpq_I} holds for the case $q=1$ as well. Lemma \ref{so(p,1) subalg lemma} below is implied in Uchida's work, but a proof is not provided. We define the following subspaces of $\mathfrak{so}(p,1)$,
\begin{equation*}
E = \text{Span} \left\{ 
    \begin{bmatrix}
         & & 1\\ & \textbf{\LARGE 0}_{p-1}& \\ 1&&
    \end{bmatrix} \right\}   
\end{equation*}
where $\textbf{0}_{p-1}$ is the $(p-1)\times (p-1)$ zero matrix, and
\begin{equation*}
    \mathfrak{t}[a:b] := \left\{ 
    \begin{bmatrix}
         0& -bu^T& 0 \\ bu& 0 &au  \\ 0&au^T&0
    \end{bmatrix} : u \in \mathbb{R}^{p-1} \right\}
\end{equation*}
where $[a:b] \in \mathbb{R}\text{P}^1$. Also, let $\mathfrak{h}$ be the Lie algebra of $H$. Then, we have:

\begin{lemma} \label{so(p,1) subalg lemma}
    Suppose $p\geq3$. Let $\mathfrak{a}$ be a proper subalgebra of $ \mathfrak{so}(p,1)$ which contains $\mathfrak{h} \simeq \mathfrak{so}(p-1)$ and suppose dim $\mathfrak{so}(p,1) -$ dim $\mathfrak{a} \leq p$. Then, $\mathfrak{a}=\mathfrak{t}[a:b]$ for some $[a:b] \in \mathbb{R}\textup{P}^1$ or $\mathfrak{a}=E \oplus \mathfrak{t}(1:\pm1)$.
\end{lemma}
\begin{proof}
Consider the following subspaces of $\mathfrak{so}(p,1)$:
    \begin{equation*}
     T_1 = \text{Span} \left\{ 
    \begin{bmatrix}
         0& -v^T& 0 \\ v& 0&0  \\ 0&0&0
    \end{bmatrix} : v \in \mathbb{R}^2 \right\}
\end{equation*}
and
\begin{equation*}
     T_2 = \text{Span} \left\{ 
    \begin{bmatrix}
         0& 0& 0 \\ 0& 0 &w  \\ 0&w^T&0
    \end{bmatrix} : w \in \mathbb{R}^2 \right\}
\end{equation*}
Then, $\displaystyle \mathfrak{so}(p,1) = \mathfrak{so}(p-1) \oplus E \oplus T_1 \oplus T_2$. Moreover, if we consider the adjoint action of $\text{SO}(p-1)$ on $\mathfrak{so}(p,1)$, all the above summands are invariant and irreducible. In addition:
\begin{align*}
    \bullet \quad & [T_2 \oplus E,T_2 \oplus E] = \mathfrak{so}(p) \\
    \bullet \quad & [T_1,E] = T_2 \\
    \bullet \quad & [T_1,T_1] = 0 
\end{align*}

 Now, $\mathfrak{a}$ is also $\text{SO}(p-1)$ invariant and $\mathfrak{a} \cap \big( T_1 \oplus T_2 \big) \neq \{ 0 \}$ because of dim $\mathfrak{so}(p,1) -$ dim $\mathfrak{a} \leq p$. Therefore $\mathfrak{a} \cap \big( T_1 \oplus T_2 \big) = \mathfrak{t}[a:b]$. Suppose $\mathfrak{a} \cap E \neq \{0\} $ which implies $ \mathfrak{a} \cap E = E$. In that case, $\mathfrak{a} \cap \big( T_1 \oplus T_2 \big) = \mathfrak{t}[a:b]$, with $a,b \neq 0$. Indeed, if $b = 0$, then $T_2 \leq \mathfrak{a}$ and hence $\mathfrak{a} = \mathfrak{so}(p,1)$, since $[T_2 \oplus E,T_2 \oplus E] = \mathfrak{so}(p)$. On the other hand, if $a=0$, then $T_1 \leq \mathfrak{a}$ hence, $ [T_1,E] = T_2 \leq \mathfrak{a} $ and as before we conclude $\mathfrak{a} = \mathfrak{so}(p,1)$. Therefore, if $\mathfrak{a} \cap E = E$, then $\mathfrak{a} = \mathfrak{so}(p-1) \oplus E \oplus \mathfrak{t}[a:b]$, with $[a:b] \in \mathbb{R}P^1$. Assume that $[a:b] \neq [1: \pm 1]$. Then
\begin{align*}
   & \begin{bmatrix}
         0& -bu^T& 0 \\ bu& 0 &au  \\ 0&au^T&0
    \end{bmatrix} \cdot \begin{bmatrix}
         & & 1\\ & \textbf{\LARGE 0}_{p-1}& \\ 1&&
    \end{bmatrix} - \begin{bmatrix}
         & & 1\\ & \textbf{\LARGE 0}_{p-1}& \\ 1&&
    \end{bmatrix} \cdot  \begin{bmatrix}
         0& -bu^T& 0 \\ bu& 0 &au  \\ 0&au^T&0
    \end{bmatrix} \\
    = &\begin{bmatrix}
        0& -au^T &0 \\ 
        au &0&bu \\
        0& bu^T&0
    \end{bmatrix} \in \mathfrak{a} \\
    \xRightarrow{\times \frac{b}{a}} & \begin{bmatrix}
        0& -bu^T &0 \\ 
        bu &0&\frac{b^2}{a}u \\
        0& \frac{b^2}{a}u^T&0
    \end{bmatrix} \in \mathfrak{a} 
    \Rightarrow  \begin{bmatrix}
        0& 0 &0 \\ 
        0 &0&\frac{a^2-b^2}{a}u \\
        0& \frac{a^2-b^2}{a}u^T&0
    \end{bmatrix} \in \mathfrak{a} \\
    \Rightarrow & T_2 \leq \mathfrak{a} \qquad \big( \text{since } [a:b] \neq [1:\pm1] \big)
\end{align*}
As before, we arrive at the contradiction that $\mathfrak{a} = \mathfrak{so}(p,1)$.
\end{proof}
 Given a $G$ action on a manifold $M$, Lemma \ref{so(p,1) subalg lemma} above gives us an analytic function $f: \mathcal{F} \rightarrow \mathbb{R}\text{P}^1$, see \cite[Remark 4]{Lentas_SOpq_I}. 

\subsection{Construction of $\text{SO}^\circ (p,1)$ actions}

\subsubsection{$(p,1)$-basic construction (I)}
The main ingredient for the constructions below will be a special type of flows on $\text{S}^1$. Here, we see $\text{S}^1$ as the unit circle in $\mathbb{R}^2$, namely $\text{S}^1 = \{ \alpha \, e_1 + \beta \, e_2 : \alpha^2+\beta^2=1 \}$.
\begin{definition} \label{basic j1 flow}
	Assume $\Phi_\theta$ is a nontrivial analytic flow on $\textup{S}^1$ and let  $\textup{J}_1$ be the reflection with respect to the $y$-axis, namely $\textup{J}_1 ( \alpha \, e_1 + \beta \, e_2) = - \alpha \, e_1 + \beta \, e_2$. We will say that $\Phi_\theta$ is a \underline{basic $\textup{J}_1$-flow} if:
	\begin{itemize}
		\item It has exactly two fixed points on $\textup{S}^1$, none of which are $\pm e_2$.
		\item $\textup{J}_1 \Phi_\theta(z) = \Phi_{-\theta} \left( \textup{J}_1 z \right)$, for $\theta \in \mathbb{R}$ and $z \in \textup{S}^1$.
		\item The Jacobian of $\Phi_\theta$ is $-\frac{2}{n}$, respectively $\frac{2}{n}$, at the attracting, respectively repelling, fixed point, where $n \in \mathbb{N}$.
	\end{itemize}
\end{definition}
Now, let $\mathcal{S}=  \left\{ \alpha e_1 + \beta e_{p+1} : \alpha^2 + \beta^2 = 1 \right\} \subseteq \text{S}^p$ and let
\begin{equation} \label{rho: S1 -> mathcal S}
    \rho_0 : \text{S}^1 \rightarrow \mathcal{S}
\end{equation}
be the analytic isomorphism defined by $\alpha \, e_1 + \beta \, e_2 \mapsto \alpha \, e_1 + \beta \, e_{p+1}$.
\begin{definition} \label{induced basic j1 flow}
    Suppose $\Phi_\theta$ is a basic $\textup{J}_1$-flow on $\textup{S}^1$. Via the isomorphism $\rho_0$, $\Phi_\theta$ induces an analytic flow $\Phi_\theta'$ on $\mathcal{S}$. We will call analytic flows on $\mathcal{S}$ arising this way, \underline{induced basic $j_1$-flows}. 
\end{definition}

 Suppose that we have a basic $\text{J}_1$-flow, which induces, via $\rho_0$ an induced basic $j_1$-flow on $\mathcal{S}$. Note that, $\mathcal{S}$ is the fixed point set of $H$ in the standard $K$ action on $\text{S}^p$. Then, analogously to the basic construction in \cite[Section 2.1]{Lentas_SOpq_I}, see also \cite{Uchida_SOpqII}, we can construct an analytic action of $G$ on $\text{S}^p$. In that action, the restricted action of $\mathcal{M}(p,1)$ on $\mathcal{S}$ is the induced basic $j_1$-flow we started with. We will refer to actions of $G$ on $\text{S}^p$ arising this way as \textit{actions from the $(p,1)$-basic construction} (I). By similar arguments as in \cite[Section 2]{Lentas_SOpq_I}
we get:
\begin{theorem}
\textup{(see also \cite[Theorem]{Uchida_SOpqII})} The analytic actions of $\textup{SO}^\circ (p,1)$ on $\textup{S}^p$ from the $(p,1)$-basic construction \textup{(I)} are analytically isomorphic if and only if the corresponding basic $\textup{J}_1$-flows are analytically isomorphic. 
\end{theorem}

\subsubsection{$(p,1)$-basic construction (II)}
 Here we construct $G$ actions on $\text{S}^{p-1} \times \text{S}^1$. However, we need to modify the definition of the flows used as the main data. This is due to the fact mentioned above, that there can now be more than two orbits with isotropy group $\text{SO}^\circ(p-1,1)$. 
\begin{definition}
    Assume $\Phi_\theta$ is a nontrivial analytic flow on $\text{S}^1$. We will say that $\Phi_\theta$ is a \underline{simple flow} if:
    \begin{itemize}
        \item The fixed points of $\Phi_\theta$ alternate between attracting and repelling.
        \item The Jacobian of $\Phi_\theta$ at a fixed point $z_i$ is $J_\Phi(z_i)= \frac{2}{n_i}$, respectively $J_\Phi(z_i)= - \frac{2}{n_i}$, if $z_i$ is an attracting, respectively repelling, fixed point, where $n_i \in \mathbb{N}$.
    \end{itemize}
\end{definition}
 Note that in the case of a simple flow, the number of fixed points must be even, see \cite{HITCHIN1991359}. Moreover, we can define a function $f_\Phi$ as in \cite[Section 2]{Lentas_SOpq_I}.  Let $p_i \in \mathcal{S}$ such that $p_i$ is between $z_i$ and $z_{i+1}$ for $1 \leq i \leq n-1$, and $p_n$ is between $z_n$ and $z_1$. Define a function $f_\Phi$ in the following way: 
\begin{itemize}
    \item $f_\Phi(p_i) := 0$, $1 \leq i \leq n$
    \item $f_\Phi(z) : = \text{tanh}(\theta)$, for $z = \Phi_\theta \left( p_i \right)$,  $\theta \in \mathbb{R}$, in the $\Phi_\theta$-orbit of a $p_i$
    \item $f_\Phi(z_i): = 1$, respectively $f_\Phi(z_i):=-1$, if $z_i$ is attracting, respectively repelling
\end{itemize}
The function $f_\Phi$ is analytic, see the proof of \cite[Lemma 2.3]{Lentas_SOpq_I}. The fixed point set of $H= \text{SO}(p-1)$ in the standard action of $K$ on $\text{S}^{p-1} \times \text{S}^1$ is
\[ \mathcal{F} = \{ \pm e_1 \} \times \text{S}^1\]
Let
\[ \mathcal{S} = \{ e_1 \} \times \text{S}^1\]
Then, identifying $\text{S}^1$ with $\mathcal{S}$, we can construct $G$ actions on $\text{S}^{p-1} \times \text{S}^1$ in a way analogous to the basic construction in \cite[Section 2.1]{Lentas_SOpq_I} using simple flows. These actions extend the standard $K$ action and the restricted $\mathcal{M}(p,1)$ action on $\mathcal{S}$ is the flow induced by the simple flow we started with, via the identification of $\text{S}^1$ with $\mathcal{S}_1$.

 We will refer to the actions on $\text{S}^{p-1} \times \text{S}^1$ obtained this way as actions \textit{from the $(p,1)$-basic construction} (II). These actions can be classified via the vector fields that the corresponding simple flows $\Phi_\theta$ induce $\text{S}^1$, similarly to \cite[Theorem 2.1]{Lentas_SOpq_I}. In particular, they can be classified by the Jacobian of $\Phi_\theta$ at the fixed points and the global invariant $\mu$, see \cite{HITCHIN1991359}.
\begin{theorem}
Let $p\geq 3$. Two $\textup{SO}^\circ(p,1)$-actions on $\textup{S}^{p-1} \times \textup{S}^1$ from the $(p,1)$-basic construction \textup{(II)} are analytically isomorphic if and only if the corresponding simple flows are analytically isomorphic. Moreover, the simple flows on $\textup{S}^1$ are classified up to analytic isomorphism by the number of fixed points, the Jacobian at the fixed points, and the global invariant $\mu$.
\end{theorem}
\begin{remark}
 When proving the analyticity of the actions from the $(p,1)$-basic construction (II), a small modification is required in the existence part of equation (33) in \cite[Appendix A]{Lentas_SOpq_I}. As we mentioned above, there can be more than two points where $f_\Phi =0$. Hence, using the notation of \cite[Appendix A]{Lentas_SOpq_I}, there are not just two points that play the role of $(e_{p+1},\epsilon_1)$ and $(-e_{p+1}, \epsilon_1)$. However, for a $z \in S_+$ with $ 0 < f(z) < 1$, there exists a $z_0 \in S_+$ such that $f_\Phi( z_0)=0$ and $z$ is in the orbit of $z_0$. Then, the proof proceeds the same.
\end{remark}
\subsection{Classification}

 Now, suppose we have a $G$ action on a manifold $M$ of dimension $p$, let $x \in M$ be a point that is not fixed by $\text{SO}(p)$. Looking at the restricted $K \simeq \text{SO}(p)$ action, let $\mathcal{O}^p$ denote the $\text{SO}(p)$-orbit of $x$. Then, dim$\mathcal{O}^p \leq p$. Using Uchida's classification of subgroups of $\text{SO}(p)$ of codimension at most $p$ ,see \cite{Uchida_clssf_rlanalytic_SLn}, and a similar analysis as in \cite[Section 3]{Lentas_SOpq_I}, we can see that necessarily $\mathcal{O}^p \simeq \text{SO}(p) \Big/ \text{SO}(p-1)$. Therefore, we can conclude that the fixed point set of $H$, $\mathcal{F}$, is nonempty and that its connected components are 1-dimensional. 

\begin{lemma}
    $G$ does not have any fixed points.
\end{lemma}
\begin{proof}
    Suppose $x \in M$ is a fixed point of $G$. Firstly, we show that $x$ is an isolated fixed point. Indeed, first of all, $x$ is also a fixed point of $\text{SO}(p)$. Note that, since $G= \text{SO}^\circ(p,1)$ is simple, its action is locally faithful, and so the $\text{SO}(p)$ action is also locally faithful. Since $\text{SO}(p)$ is compact, its action is linearisable at $x$, and by choosing an $\text{SO}(p)$-invariant metric for $M$, we get a representation
    \begin{equation*}
        \rho_1 : \text{SO}(p) \rightarrow \text{O}(p)
    \end{equation*}
    This representation must be either trivial or irreducible. Because if there was a $k$-dimensional, $1 \leq k < p$, $\rho_1$-invariant subspace of $\mathbb{R}^p$, $1 \leq k < p$, then (changing basis if necessary) we'd get a homomorphism 
    \begin{equation*}
        \text{SO}(p) \hookrightarrow \text{O}(k) \times \text{O}(p-k)
    \end{equation*}
    Since the action is locally faithful, Ker$\rho_1$ must be $0$-dimensional. Comparing dimension we get $\displaystyle \frac{p(p-1)}{2} \leq \frac{k(k-1)}{2} + \frac{(p-k)(p-k-1)}{2}$, which yields a contradiction. Now, if $\rho_1$ is trivial, that means the action of $\text{SO}(p)$ is trivial, and consequently the action of $G$ is trivial. So, $\rho_1$ is irreducible and $x$ is an isolated fixed point of $\text{SO}(p)$ and hence of $G$.
    
     Consider the isotropy representation
    \begin{equation*}
        \rho : \mathfrak{so}(p,1) \rightarrow \mathfrak{gl}_p(\mathbb{R})
    \end{equation*}
    Since $\mathfrak{so}(p,1)$ is simple, and the representation is not trivial, it has to be irreducible. Complexifying and using that $\displaystyle \mathfrak{so}(p,1) \cong \mathfrak{so}(p+1)(\mathbb{C})$, we get an irreducible representation 
    \begin{equation*}
        \rho_1^{\mathbb{C}}: \mathfrak{so}(p+1)(\mathbb{C}) \rightarrow \mathfrak{gl}_{p}(\mathbb{C}) 
    \end{equation*}
    See \cite{onishchik2004lectures}.
    Using the Weyl character formula,  \cite[p.408-410]{FultonHarris}, we can see that the dimension of an irreducible representation is $\displaystyle \binom{p+1}{k}$ if $p+1$ is odd, and $\displaystyle \binom{p+1}{k}$ or $\displaystyle \frac{1}{2} \binom{p+1}{k}$ if $p+1$ is even. Comparing $p$ to these binomial coefficients, we get a contradiction. 
\end{proof}

In particular, $G$ does not have a fixed point on $\mathcal{F}$ and hence by Lemma \ref{so(p,1) subalg lemma} we can define an analytic function 
\[ f: \mathcal{F} \rightarrow \mathbb{R} \text{P}^1\]
and by the action of $\mathcal{M}(p,1)$ on $\mathcal{F}$ we get a flow
\[ \Phi : \mathbb{R} \times \mathcal{F} \rightarrow \mathcal{F}\]

 Then, as in \cite[Section 4]{Lentas_SOpq_I} we can show the following:



\begin{theorem}
    Suppose $\textup{SO}^\circ(p,1)$ acts analytically on a closed connected manifold $M$ of dimension $p$. Then, if $\textup{SO}(p)$ has a fixed point, $M$ is equivariantly covered by $\textup{S}^p$, where the $\textup{SO}^\circ(p,1)$ action on $\textup{S}^p$ is one from the $(p,1)$-basic construction \textup{(I)}. If $\textup{SO}(p)$ does not have a fixed point, then
    \begin{itemize}
        \item If there exists a nullcone orbit, then $M$ is equivariantly covered by $\textup{SO}^\circ(p,1) \times_P \textup{S}^1$ with the standard left $\textup{SO}^\circ(p,1)$ action, where where $P \leq \textup{SO}^\circ(p,1)$ is a maximal parabolic subgroup isomorphic to the stabiliser of an isotropic line in the standard representation of $\textup{SO}^\circ(p,1)$ on $\mathbb{R}^{p+1}$. If $P=M_P A_P N_P$ is a Langlands decomposition of $P$, $P$ acts on $\textup{S}^1$ by a flow via $A_P$.
        \item If there does not exist a nullcone orbit, then $M$ is equivariantly covered by $\textup{S}^{p-1} \times \textup{S}^1$, where the $\textup{SO}^\circ(p,1)$ action on $\textup{S}^{p-1} \times \textup{S}^1$ is one from $(p,1)$-basic construction \textup{(II)}.
    \end{itemize}
\end{theorem}
\section{Classification of analytic $\text{SO}^o(p,2)$, $p\geq3$, actions} \label{SO(p,2) axns}

Here, we consider the case $q =2$ and $p \geq 3$, namely the analytic actions of $G= \text{SO}(p,2)$ on manifolds of dimension $p+1  = p+q-1$. Uchida in \cite{Uchida_SOpqII} studied these kind of actions on the sphere $\text{S}^{p+1}$. Modifying Uchida's methods, we get actions of $G$ on $\text{S}^p \times \text{S}^1$ and $\text{S}^{p-1} \times \text{S}^2$ where $K= \text{SO}(p) \times \text{SO}(2)$ acts in a standard way, see below. Then we see that an action of $G$ on a closed manifold of dimension $p+1$ is equivariantly covered by one of those above.

Consider the following actions of $K$ on $\text{S}^p \times \text{S}^1$ and $\text{S}^{p-1} \times \text{S}^2$. For $\text{S}^p \times \text{S}^1$, we embed $\text{SO}(p)$ in $\text{SO}(p+1)$ as top left block matrices.
Then, via $\text{SO}(p+1)$, $\text{SO}(p)$ acts on $\text{S}^p$. The action of $\text{SO}(2)$ on $\text{S}^1$ is the standard one and that gives the action of $K$ on $\text{S}^p \times \text{S}^1$. Similarly, we get the action of $K$ on $\text{S}^{p-1} \times \text{S}^2$. We will refer to these actions as \textit{standard} actions of $K$ on these spaces.

The main difference with the $q \geq 3$ case is that the fixed point set of $H   = \text{SO}(p-1) \times \text{SO}(1)  \simeq \text{SO}(p-1)$ is no longer 1-dimensional, but 2-dimensional. Abusing the notation we will just write $H = \text{SO}(p-1)$. We will need some more terminology. Following \cite{Uchida_SOpqII}, let $N$ be the subgroup of $G$ with Lie algebra \[ \mathfrak{n} =  \left\{ \begin{bmatrix}
    0 & & \alpha & \beta \\
    & \textbf{\LARGE 0}_{p-1} && \\
    \alpha && 0& \gamma \\
    \beta && - \gamma & 0
\end{bmatrix} : \alpha, \beta, \gamma \in \mathbb{R} \right\} \leq \mathfrak{so}(p,2)\]
where $\textbf{\LARGE 0}_{p-1}$ is the zero square matrix of dimension $p-1$ and any omitted value is assumed to be equal to $0$. Note that, $\text{SO}(2) \leq N$, where by $\text{SO}(2)$ we mean the subgroup $ \{I_p\} \times\text{SO}(2)$. Moreover, let 
$$N^+ = N \, \bigcup \, j_1N$$
Note that $j_1j_2$ commutes with $N$.

 Let $\mathcal{F}$ be the fixed point set of $H$ and $f: \mathcal{F} \rightarrow \mathbb{R} \text{P}^2$. For $Y \in \mathcal{F}$ with $f(Y) = [a :b:c]$, let
 \begin{equation} \label{U_Y RP2 def eq}
      U^\circ(Y) =  H_{[a:b:c]}^\circ
 \end{equation}
where $H_{[a:b:c]}$ is the stabiliser of the point $a e_1 + b e_{p+1} + c e_{p+2}$ in the standard representation of $G$ on $\mathbb{R}^{p+2}$. The following equality also holds
\[ G = KNU^\circ(Y)\]
See \cite{Uchida_SOpqII}.
Moreover, the subgroup
\[ \mathcal{M}(p,2) = \left\{ \begin{bmatrix}
    \text{cosh}(\theta) & & \text{sinh}(\theta) & \\
    & I_{p-1} && \\
    \text{sinh}(\theta) & & \text{cosh}(\theta) & \\
    &&& 1
\end{bmatrix} \right\}\]
is in $N$ and we have the following equality
\begin{equation} \label{N decomp}
    N = \text{SO}(2) \, \mathcal{M}(p,2) \, (U^\circ(Y) \cap N)
\end{equation}
for any $Y \in \mathcal{F}$, as is established in \cite{Uchida_SOpqII}. Note that (\ref{N decomp}) is equivalent to
\begin{equation} \label{N decomp with H(a b c)}
    N = \text{SO}(2) \, \mathcal{M}(p,2) \, ( H_{[a:b:c]}^\circ \cap N)
\end{equation}
Note that $ N \simeq \text{SO}^\circ(1,2)$,  the obvious isomorphism $\varphi : N \rightarrow \text{SO}^\circ(1,2)$ being
\[  \begin{bmatrix}
    x_1 & & x_2 & x_3 \\
    & O_{p-1} && \\
    x_4 && x_5& x_6 \\
    x_7 && x_8 & x_9
\end{bmatrix}  \mapsto  \begin{bmatrix}
    x_1 &  x_2 & x_3 \\
    x_4 &x_5& x_6 \\
    x_7 & x_8 & x_9
\end{bmatrix} \]
It is well known that $\text{SO}^\circ(1,2)$ is isomorphic to  $\text{PSL}_2(\mathbb{R})$ and hence, $N$ is also isomorphic to $\text{PSL}_2(\mathbb{R})$. Let $\mathcal{M}(1,2)$ be the subgroup of $\text{SO}^\circ(1,2)$
\begin{equation} \label{M(1,2) def}
\mathcal{M}(1,2) = \left\{ \begin{bmatrix}
    \text{cosh}(\theta) &  \text{sinh}(\theta) & \\
    \text{sinh}(\theta) &  \text{cosh}(\theta) & \\
    && 1
\end{bmatrix} : \theta \in \mathbb{R}\right\}.
\end{equation}
Then, $\varphi$ gives an isomorphism between $\mathcal{M}(p,2)$ and $\mathcal{M}(1,2)$. We will denote an element of either group by $m(\theta)$. It will be clear from the context every time what group it belongs to.

Moreover, via $\varphi$, we get an equation analogous to (\ref{N decomp with H(a b c)}) for $\text{SO}^\circ(1,2)$. Let $\widetilde{H}_{[a:b:c]}$ be the stabiliser of the point $a e_1 + b e_{2} + c e_{3}$ in the standard representation of  $\text{SO}^\circ(1,2)$ on $\mathbb{R}^{3}$. Then, we have
\begin{equation} \label{SO(1,2) decomp with H(abc)}
    \text{SO}^\circ(1,2) = \text{SO}(2) \, \mathcal{M}(1,2) \,\widetilde{H}_{[a:b:c]}^\circ 
\end{equation}
 Now, let $j_1' = \begin{bmatrix}
    -1&& \\ &0& \\ &&0
\end{bmatrix}$ and consider the group
\begin{equation} \label{SO(1,2)+ def}
    \text{SO}^\circ(1,2)^+ = \text{SO}^\circ(1,2) \cup j_1' \,\text{SO}^\circ(1,2)
\end{equation}
Defining $\varphi(j_1) = j_1'$, it is immediate that $\varphi$ extends to an isomorphism between $N^+$ and $\text{SO}^\circ(1,2)^+$
\begin{equation} \label{N^+ SO(1,2)+ isom}
    \varphi : N^+ \xrightarrow{\simeq} \text{SO}^\circ(1,2)^+
\end{equation}

In the case of $\text{SO}^\circ(p,2)$ actions there exists a lemma similar to \cite[Lemma 2.2]{Lentas_SOpq_I}. We denote by $\mathfrak{h}_{[a:b:c]}$ the Lie algebra of the group $H_{[a:b:c]}$ defined earlier.
\begin{lemma}{\cite[Lemma 1.2]{Uchida_SOpqII}} \label{SOp,2 lie alg lemma}
    Suppose $p \geq 3$ and let $\mathfrak{a}$ be a proper Lie subalgebra of $\mathfrak{so}(p,2)$ which contains $\mathfrak{h} \simeq \mathfrak{so}(p-1)$. If 
    \[ \textup{dim} \,  \mathfrak{so}(p,2) - \textup{ dim} \, \mathfrak{a} \leq p+1\]
    then
    \[ \mathfrak{a}= \mathfrak{h}_{[a:b:c]}\]
    for some $(a,b,c) \neq (0,0,0)$, or 
    \[ \mathfrak{a}= \mathfrak{h}_{[a:b:c]} \oplus R^1\]
    for some $(a,b,c) \neq (0,0,0)$ such that $a^2=b^2+c^2$, where the space $R^1$ is generated by the matrix $b \left( E_{1,p+1} + E_{p+1,1} \right) + c \left( E_{1,p+2} + E_{p+2,1} \right) $.

\begin{remark} \label{f RP2 def remark}
Similarly to \cite[Remark 4]{Lentas_SOpq_I}, Lemma \ref{SOp,2 lie alg lemma} above allows us to define an analytic function $f:\mathcal{F} \rightarrow \mathbb{R}\text{P}^2$, where $\mathcal{F}$ is the fixed point set of $H$, by the property that for $Y \in \mathcal{F}$, $f(Y) \in \mathbb{R}\text{P}^2$ is the unique point such that $\mathfrak{h}_{f(Y)} \leq \mathfrak{g}_Y$, where $\mathfrak{g}_Y$ is the isotropy algebra of $Y$ with respect to the $G$ action.
 \end{remark}

\end{lemma}

\subsection{Extracting the main data}
In this section, we begin with a $G$ action on $\text{S}^p \times \text{S}^1$, $\text{S}^{p-1} \times \text{S}^2 $, or $\text{S}^{p+1}$ extending the standard $K$ action in each case. We will show that to such $G$ actions correspond an $\text{SO}^\circ(1,2)^+$ action on $\mathbb{T}^2$ in the first case, an $\text{SO}^\circ (1,2)$ action on $\text{S}^2$ in the second case, and an $\text{SO}^\circ(1,2)^+$ action on $\text{S}^2$ in the third case. In Section \ref{(p,2) case constr} we will show that these actions are all we need in order to construct and classify the $G$ actions that extend the standard $K$ action in each case, see in particular Theorems \ref{bsc constr (I) thm}, \ref{bsc constr (II) thm}, and \ref{bsc constr (III) thm}.

\begin{remark}[See \cite{Uchida_SOpqII}] \label{SO(p,2) axn -> (phi,f)}
Let $M$ be $\textup{S}^{p+1}$, $\textup{S}^p \times \textup{S}^1$, or $\textup{S}^{p-1} \times \textup{S}^2$. Suppose $G$ acts analytically on $M$ in a way that extends the standard $K$ action in each case. Consider the fixed point set, $\mathcal{F}$, of $H$. Let $\phi$ be the restricted action of $N^+$ on $\mathcal{F}$ and $f:\mathcal{F} \rightarrow \mathbb{R}\text{P}^2$ as in Remark \ref{f RP2 def remark}. Then, the pair $(\phi,f)$ satisfies:
\begin{itemize}
    \item[(B1)] $\phi$ is an analytic action of $N^+$ on $\mathcal{F}$, such that the restricted action of $N^+ \bigcap K$ is the restriction of the standard action of $K$.
    \item[(B2)] $f$ is an analytic function from $\mathcal{F} \rightarrow \mathbb{R}\text{P}^2$, which is $N^+$-equivariant and satisfies: $N^+ \bigcap U^\circ(Y) \subset N^+_Y$
\end{itemize}
Moreover, let $\mathcal{R} = \left\{ [a:b:c] \in \mathbb{R}\text{P}^2 : c=0 \right\} \subseteq \mathbb{R}\text{P}^2$ and $\rho : \mathcal{R} \rightarrow \mathbb{R}\text{P}^1$ the obvious isomorphism. Consider $\mathcal{S}=f^{-1}(\mathcal{R})$. Then, $\mathcal{S}$ is a one dimensional submanifold of $\mathcal{F}$ which is $j_1$ and $j_2$ invariant and it is transverse to each $\text{SO}(2)$-orbit, see \cite{Uchida_SOpqII}.  Let $\Phi_\theta$ be the flow that $\mathcal{M}(p,2)$ induces on $\mathcal{S}$ and $f_0 = \rho \circ f\big| _\mathcal{S}$. Then, $\Phi_\theta$ is an analytic flow on $\mathcal{S}$ , $f_0:\mathcal{S} \rightarrow \mathbb{R}\text{P}^1$ is analytic, and they satisfy 
\begin{itemize}
    \item[(i)] $j_i \star \Phi_\theta (z) = \Phi_{- \theta} (j_i \star z)$ \qquad $(i=1,2)$
    \item[(ii)] $f_0(z) = [a:b:0] \Rightarrow f_0(j_i \star z) = [a : -b:0]$ \qquad $(i=1,2)$
    \item[(iii)] $f_0(z) = [a:b:0] \Rightarrow f_0 \left( \Phi_\theta (z) \right) = \left[ a \, \text{cosh}(\theta) + b \, \text{sinh}(\theta) : a \, \text{sinh}(\theta) + b \, \text{cosh}(\theta) :0\right] $ 
    \item[(iv)] $f_0(z) = [0:1:0] \Leftrightarrow j_1 \star z=z \quad \text{and} \quad f_0(z) = [1:0:0] \Leftrightarrow j_2 \star z = z$
\end{itemize}
These are analogous to the Uchida conditions from \cite[Section 2]{Lentas_SOpq_I}, see \cite[Remark 7]{Lentas_SOpq_I}.
\end{remark}

\subsubsection{Case 1: $G$ actions on $\text{S}^p\times \text{S}^1$} \label{SO(p,2) axns on Sp x S1}
Here we will see that a $G$ action on $\text{S}^p \times \text{S}^1$ extending the standard $K$ action correspond to a special type of $\text{SO}^\circ(1,2)^+$, see Definition \ref{SO(1,2)+ def}, actions on $\mathbb{T}^2$ defined below.

\begin{definition} \label{T2-basic axns}
    Assume that $\phi$ is an $\text{SO}^\circ(1,2)^+$ action on $\mathbb{T}^2=\text{S}^1 \times \text{S}^1$. We will say that $\phi$ is a \underline{$\mathbb{T}^2$-basic action} if
    \begin{itemize}
        \item $j_1'$ acts as the reflection with respect to the $y$-axis on the first factor and trivially on the second factor of $\mathbb{T}^2$.
        \item $\text{SO}(2)$ acts trivially on the first factor and in the usual way on the second factor.
        \item $\phi$ has two closed circle orbits none of which are the circles fixed by $j_1'$.
    \end{itemize}
\end{definition}

Suppose that $G$ acts analytically on $\text{S}^p \times \text{S}^1$ and the restricted $K$ action is the standard one. Here, the fixed point set, $\mathcal{F}$, of $H$ is
$$\mathcal{F} = \left\{ \alpha e_1 + \beta e_{p+1} : \alpha^2 + \beta^2=1 \right\} \times \text{S}^1 \simeq \mathbb{T}^2$$ 
Lemma \ref{SOp,2 lie alg lemma} implies the existence of a function $f: \mathcal{F} \rightarrow \mathbb{R}\text{P}^2$. Denote the restricted action of the subgroup $N^+$ on $\mathcal{F}$ by $\phi_1$. Then, the pair $(\phi_1 , f)$ satisfies properties (B1)-(B2). Recall that $N^+$ is isomorphic to $\text{SO}^\circ(1,2)^+$, see the beginning of this section. Therefore, identifying $\mathcal{F}$ with $\mathbb{T}^2$ in the obvious way and using the isomorphism between $N^+$ and $\text{SO}^\circ(1,2)^+$, we get an action $\phi$ of $\text{SO}^\circ(1,2)^+$ on $\mathbb{T}^2$. Note that $N^+ \cap K = \text{SO}(2) \cup j_1 \text{SO}(2)$, hence $\phi$ is a $\mathbb{T}^2$-basic action. It is immediate that starting with two analytically isomorphic $G$ actions will result in analytically isomorphic $\mathbb{T}^2$-basic action of   $\text{SO}^\circ(1,2)^+$ on $\mathbb{T}^2$. \\

\begin{remark} \label{Schneider T2}
Schneider, in \cite{SchneiderSL2axns}, classified $\text{SL}_2(\mathbb{R})$ actions on $\mathbb{T}^2$. Fix a Borel subgroup $\Sigma$ of $\text{SL}_2(\mathbb{R})$ and an infinitesimal generator, say $\alpha$, of the diagonal subgroup. Then, the $\text{SL}_2(\mathbb{R})$ actions come down to a 1-dimensional, closed submanifold $W$ of $\mathbb{T}^2$ and 
an action of $\Sigma$ on $W$. Furthermore, 
\[ \mathbb{T}^2 \simeq \text{SL}_2(\mathbb{R}) \times_\Sigma W\]
and the action is the left $\text{SL}_2(\mathbb{R})$ action on $\text{SL}_2(\mathbb{R}) \times_\Sigma W$ (See \cite[Theorem 1]{SchneiderSL2axns}). Two such actions of $\text{SL}_2(\mathbb{R})$ are isomorphic, if the vector fields on $W$ generated by $\alpha \in \mathfrak{sl}_2(\mathbb{R})$ are isomorphic.
\end{remark}
Here, $\text{SL}_2(\mathbb{R}) \times_\Sigma W$ is the quotient space of $\text{SL}_2(\mathbb{R}) \times W$ modulo the equivalence relation
\[ (\mu , w) \sim (\mu \cdot \sigma , \sigma^{-1} \star w )  \]
where $\mu \in \text{SL}_2(\mathbb{R}), \sigma \in \Sigma$ and $w \in W$.

Since $\text{SO}^\circ(1,2)$ is isomorphic to $\text{PSL}_2(\mathbb{R})$, by Schneider's result the analytic actions of $\text{SO}^\circ(1,2)$ on $\mathbb{T}^2$ are also classified. It is also immediate that two analytic $\text{SO}^\circ(1,2)^+$ actions on $\mathbb{T}^2$ are isomorphic if and only if the restricted $\text{SO}^\circ(1,2)$ actions are isomorphic via an isomorphism that commutes with $j_1'$.

\subsubsection{Case 2: $G$ actions on $\text{S}^{p-1}\times \text{S}^2$}

 Suppose we have a $G$ action on $\text{S}^{p-1} \times \text{S}^2$ extending the standard $K$ action. We will see that such an action corresponds to a special type of $\text{SO}^\circ(1,2)$ actions on $\text{S}^2$.

\begin{definition} \label{S2 basic axn no j1}
    Assume $\phi$ is an analytic $\text{SO}^\circ(1,2)$ action on $\text{S}^2$. We will say that $\phi$ is an \underline{$ \text{S}^2$-basic action} if: 
    \begin{itemize}
        \item the action of $\text{SO}(2)$ on $\text{S}^2$ is rotation with respect to the $z$-axis.
        \item $\phi$ has one closed circle orbit and no fixed points.
    \end{itemize}
\end{definition}

For a $G$ action on $\text{S}^{p-1} \times \text{S}^2$ like above, the fixed point set, $\mathcal{F}$, of $H$ is
\[ \mathcal{F} = \left\{ \pm e_1 \right\} \times \text{S}^2\]
Let 
\[ \mathcal{F}_1 = \left\{  e_1 \right\} \times \text{S}^2\]
and $\mathcal{F}_2$ be the other connected component of $\mathcal{F}$. We note that the action of $N$ on $\mathcal{F}_2$ is determined by the action of $N$ on $\mathcal{F}_1$ by the equation
\[ n \star (j_1 \star v) = j_1 \left( j_2nj_2 \right) \star v\]
for $n \in N$, $v \in \mathcal{F}_1$. By Lemma \ref{SOp,2 lie alg lemma}, there exists a function $f: \mathcal{F} \rightarrow \mathbb{R}\text{P}^2$ and if we denote the restricted action of the subgroup $N^+$ on $\mathcal{F}$ by $\phi_1$, then the pair $(\phi_1 , f)$ satisfy properties (B1)-(B2). Let $\mathcal{S}=f^{-1}(\mathcal{R})$, where $\mathcal{R} = \left\{ [a:b:c] \in \mathbb{R}\text{P}^2 : c=0 \right\} \subseteq \mathbb{R}\text{P}^2$. As in Remark \ref{SO(p,2) axn -> (phi,f)}, it can be shown that $\mathcal{S}$ is a closed, one dimensional submanifold of $\mathcal{F}$, which is transverse to the $\text{SO}(2)$ orbits. By the last condition, the points $(e_1, \pm \epsilon_1)$ are in $\mathcal{S}$, since they are fixed by $\text{SO}(2)$. Therefore there exists a connected component, $\mathcal{S}_1$, of $\mathcal{S}$ such that $\mathcal{S}_1 \subseteq \mathcal{F}_1$. Hence, $\mathcal{S}_1$ is a one dimensional closed submanifold of $\mathcal{F}_1$, transverse to the $\text{SO}(2)$ orbits. Now, $f \neq [0:1:0]$ since there are no $\text{SO}(p)$ fixed points in the kind of $G$ actions considered in this section. Hence, the flow, $\Phi_\theta$ that the action of $\mathcal{M}(p,2)$ induces on $\mathcal{S}$ has two fixed points, otherwise $j_2$ would have more than two fixed points, which is impossible for an involution on a one dimensional closed submanifold. Since the points $(e_1, \pm \epsilon_1)$ are in $ \mathcal{S}$, $\mathcal{S}$ is the union of the $\Phi_\theta$ orbits of these points, and the fixed points. Note that $j_2$ must map the fixed points of $\Phi_\theta$ to each other, since $j_2 \star \Phi_\theta(z)=\Phi_{- \theta}(j_2 \star z)$ for $z\in \mathcal{S}_1$. Hence, on $\mathcal{F}_1$, $\text{SO}^\circ(1,2)$ has one closed orbit, namely the orbit of the fixed points of $\Phi_\theta$ on $\mathcal{S}_1$.\\

 Now, let $\phi_2$ be the restricted $N$ action on $\mathcal{F}_1$ and identify $\mathcal{F}_1$ with $\text{S}^2$ in the obvious way. Recall that $N$ is isomorphic to $\text{SO}^\circ(1,2)$ and that $N \cap K = \text{SO}(2)$. Then, it is immediate that we get an action $\phi$ of $\text{SO}^\circ(1,2)$ on $\text{S}^2$ which is an $ \text{S}^2$-basic action.

\begin{remark} \label{Schneider S2}
Recall that $\text{SO}^\circ(1,2)$ is isomorphic to $\text{PSL}_2(\mathbb{R})$ and hence covered by $\text{SL}_2(\mathbb{R})$. The analytic $\text{SL}_2(\mathbb{R})$ actions on $\text{S}^2$ are classified in \cite{SchneiderSL2axns}, by the number of closed orbits and the \textit{normal invariant} at each closed orbit, see \cite[Theorem 6]{SchneiderSL2axns}. Therefore, also the $\text{SO}^\circ(1,2)$ on $\text{S}^2$ are classified. The normal invariant would be equal to the Jacobian of $\Phi_\theta$ at a fixed point on $\mathcal{S}_1$, in the setting above.
\end{remark} 

\subsubsection{Case 3: $G$ actions on $\text{S}^{p+1}$}

 Finally, we consider $G$ actions on $\text{S}^{p+1}$ extending the standard $K$ action, namely the actions studied in \cite{Uchida_SOpqII}. Here, the fixed point set of $H$ is 
\[ \mathcal{F} = \{ \alpha \, e_1 + \beta \, e_{p+1} + \gamma \, e_{p+1} \, : \, \alpha^2 + \beta^2 + \gamma^2 =1 \} \simeq \text{S}^2\]
Recall $\text{SO}^\circ(1,2)^+$ from Definition \ref{SO(1,2)+ def}.
\begin{definition} \label{basic S2 axn + j1}
    Let $\phi$ be an analytic $\text{SO}^\circ(1,2)^+$ action on $\text{S}^2$. We will say that $\phi$ is an \underline{$(\text{S}^2)^+$-basic action} if: 
    \begin{itemize}
        \item $j_1'$ acts on $\text{S}^2$ as the reflection with respect to the $xy$-plane.
        \item  $\text{SO}(2)$ acts on $\text{S}^2$ as rotations with respect to the $z$-axis.
        \item $\phi$ has two closed circle orbit and no fixed points.
    \end{itemize}
\end{definition}

Suppose we have a $G$ action on $\text{S}^{p+1}$ extending the standard $K$ action. Similar to Case 2, via the isomorphism $N^+ \simeq \text{SO}^\circ(1,2)^+$ and the obvious identification of $\mathcal{F}$ with $\text{S}^2$, it can be shown that to the $G$ action corresponds an $(\text{S}^2)^+$-basic action $\phi$ of $\text{SO}^\circ(1,2)^+$ on $\text{S}^2$.
 
\begin{remark} \label{schneider S2+}
Two $(\text{S}^2)^+$-basic actions are analytically isomorphic if and only if the restricted $\text{SO}^\circ(1,2)$ actions are analytically isomorphic via an isomorphism that commutes with $j_1'$. We also note again that $\text{SO}^\circ(1,2)$ is isomorphic to $\text{PSL}_2(\mathbb{R})$ and hence covered by $\text{SL}_2(\mathbb{R})$ and the analytic $\text{SL}_2(\mathbb{R})$ actions on $\text{S}^2$ are classified in \cite{SchneiderSL2axns}.
 \end{remark}

\subsection{Construction of $\text{SO}^\circ (p,2)$ actions} \label{(p,2) case constr}

 In this section we construct $\text{SO}^\circ(p,2)$ actions on $\text{S}^{p+1}$, $\text{S}^p \times \text{S}^1$, and $\text{S}^{p-1}\times \text{S}^2$ that extend the standard $K$ action in each case. We are going to use the following two lemmas. Uchida shows the following results in the case of smooth actions of $G$ on the sphere $\text{S}^{p+1}$ whose restricted $K$ action is the standard orthogonal one. However, the arguments and techniques of \cite[Section 2]{Lentas_SOpq_I} can be applied to prove them in the analytic setting, see also \cite{Lentas_thesis}.
\begin{lemma} \label{(phi , f) -> action lemma}
Let $M$ be $\textup{S}^{p+1}$, $\textup{S}^p \times \textup{S}^1$, or $\textup{S}^{p-1} \times \textup{S}^2$. Assume $K$ acts on $M$ in the standard way in each case and let $\mathcal{F}$ be the fixed point set of $H$. Suppose:
\begin{itemize}
    \item[\textup{(B1)}] $\phi$ is an analytic action of $N^+$ on $\mathcal{F}$, such that the restricted action of $N^+ \bigcap K$ is the restriction of the standard action of $K$.
    \item[\textup{(B2)}] $f$ is a function from $\mathcal{F} \rightarrow \mathbb{R}\text{P}^2$, which is $N^+$-equivariant and satisfies: $N^+ \bigcap U^\circ(Y) \subset N^+_Y$
\end{itemize}
Here, $\mathbb{R} \text{P}^2$ is seen as a $N^+$-space, via the identification $N \simeq \text{SO}^\circ(1,2)$ and by defining the action of $j_1$ by $j_1 [a:b:c] = [-a:b:c]$. For the definition of $U^\circ(Y)$ see (\ref{U_Y RP2 def eq}). Then, there exists an analytic $G$ action on $M$ extending the standard $K$ action such that the restricted $N^+$ action on $\mathcal{F}$ is $\phi$ and $f$ is exactly the function from Remark \ref{f RP2 def remark} for this $G$ action.
\end{lemma}
\begin{lemma} \label{S f0 Phi -> (phi,f)}
Let $M$ be $\textup{S}^{p+1}$, $\textup{S}^p \times \textup{S}^1$, or $\textup{S}^{p-1} \times \textup{S}^2$. Let $\mathcal{F}$ be the fixed point set of $H$ and suppose $\mathcal{S}$ is a connected, one dimensional submanifold of $\mathcal{F}$ which is $j_1$ and $j_2$ invariant in the case $M =\textup{S}^{p+1}$, $j_1$ invariant in the case $M=\textup{S}^p \times \textup{S}^1$, and $j_2$ invariant in the case $M=\textup{S}^{p-1} \times \textup{S}^2$. Assume further that in any case $\mathcal{S}$ is transverse to each $\text{SO}(2)$-orbit. Consider $\mathcal{R} = \left\{ [a:b:c] \in \mathbb{R}\text{P}^2 : c=0 \right\}$ and let $\rho:\mathcal{R} \rightarrow \mathbb{R}\text{P}^1$ be the obvious isomorphism. Suppose $\Phi_\theta$ is an analytic flow on $\mathcal{S}$ and $f_0:\mathcal{S} \rightarrow \mathbb{R}\text{P}^1$ an analytic function such that they satisfy 
\begin{itemize}
    \item[\textup{(i)}] $j_i \star \Phi_\theta (z) = \Phi_{- \theta} (j_i \star z)$ \qquad $(i=1,2)$
    \item[\textup{(ii)}] $f_0(z) = [a:b] \Rightarrow f_0(j_i \star z) = [a : -b]$ \qquad $(i=1,2)$
    \item[\textup{(iii)}] $f_0(z) = [a:b] \Rightarrow f_0 \left( \Phi_\theta (z) \right) = \left[ a \, \text{cosh}(\theta) + b \, \text{sinh}(\theta) : a \, \text{sinh}(\theta) + b \, \text{cosh}(\theta) \right] $ 
    \item[\textup{(iv)}] $f_0(z) = [0:1] \Leftrightarrow j_1 \star z=z \quad \text{and} \quad f_0(z) = [1:0] \Leftrightarrow j_2 \star z = z$
\end{itemize}
Define an action $\phi$ of $N^+$ on $\mathcal{F}$ in the following way: let $n \in N^+$ and $Z \in \mathcal{F}$. There exist $k_1 \in \text{SO}(2)$ and $Y \in \mathcal{S}$ such that $k_1 \star Y = Z$. We then have $nk_1=j_1^\sigma n'$, where $\sigma = 0$ or $1$ and $n' \in N$. Write
\[n =  k \, m(\theta) \, u\]
with $k \in \text{SO}(2)$, $\theta \in \mathbb{R}$, and $ u \in H^\circ_{[1:1:0]} \cap N$, according to \textup{(\ref{SO(1,2) decomp with H(abc)})}. Then define
\[n \star Z := k \star \Phi_\theta(Y)\]
Finally, define a function $f: \mathcal{F} \rightarrow \mathbb{R}\text{P}^2$ by 
\begin{equation} \label{f0 eq}
      f(k\star z) := k \star \left( \rho^{-1} \circ f_0(z) \right)
\end{equation}
for $k \in \text{SO}(2)$ and $z \in \mathcal{S}$. Then, $\phi$ is an analytic action of $N^+$ on $\mathcal{F}$, $f$ is analytic, and $\phi$ and $f$ satisfy conditions \textup{(B1)-(B2)} in Remark \ref{SO(p,2) axn -> (phi,f)}.
\end{lemma}
\begin{remark} \label{G axns <-> (phi,f)}
    Let $M$ be as in the lemmas above and suppose we have a pair $(\phi, f)$ as in Lemma \ref{(phi , f) -> action lemma}. By that lemma we get a $G$ action on $M$ extending the standard $K$ action. By Remark \ref{SO(p,2) axn -> (phi,f)}, we get a pair $(\phi',f')$. Then, $\phi'=\phi$ and $f'=f$. And conversely, starting from a $G$ action on $M$, we get a pair $(\phi,f)$ satisfying (B1)-(B2), and the $G$ action coming from Lemma \ref{(phi , f) -> action lemma} is the one we started with. The same property holds for $G$ actions on $M$ extedning the standard $K$ action and triplets $(\mathcal{S},\Phi_\theta,f_0)$ as in Lemma \ref{S f0 Phi -> (phi,f)}.
\end{remark}

In the constructions below we will have an action of $N^+$ on $\mathcal{F}$ and our goal will be to apply Lemma \ref{(phi , f) -> action lemma} in order to get an action of $G$. To use that lemma we will need to construct the function $f$. To that end, we will find in every case a submanifold $\mathcal{S}$ and construct a function $f_0$, as in Lemma \ref{S f0 Phi -> (phi,f)} so that we can apply Lemma \ref{S f0 Phi -> (phi,f)}.

\subsubsection{$(p,2)$-basic construction (I)} \label{(p,2)-basic const I}

Here, we show how from a $\mathbb{T}^2$-basic action of $\text{SO}^\circ(1,2)^+$ on $\mathbb{T}^2 = \text{S}^1 \times \text{S}^1$ we can construct an action of $G$ on $\text{S}^p \times \text{S}^1$ extending the standard $K$ action. Note that, for the standard $K$ action the fixed point set of $H$ is
\[\mathcal{F} = \left\{ \alpha e_1 + \beta e_{p+1} : \alpha^2 + \beta^2=1 \right\} \times \text{S}^1\]

 Let $\phi_1$ be a $\mathbb{T}^2$-basic action of $\text{SO}^\circ(1,2)^+$ on $\mathbb{T}^2$, see Definition \ref{T2-basic axns}. Consider the function $\tau:\text{S}^1 \rightarrow \text{S}^p $ defined by
\[(x,y)\mapsto x \, e_1 + y \, e_{p+1} \]
Then, 
\begin{equation} \label{T2 and F isomor}
    \psi = \tau \times \text{Id}: \mathbb{T}^2 \rightarrow \mathcal{F}
\end{equation}
is an analytic isomorphism. Via the isomorphism $\psi$, we get an action $\phi_2$ of $\text{SO}^\circ(1,2)^+$ on $\mathcal{F}$ defined by the equation
\[ \phi_2(n,\psi(p)) = \psi \left( \phi_1(n,p) \right)\]
for $p\in \mathbb{T}^2$ and $n \in\text{SO}^\circ(1,2)^+$. Note that, since $\phi_1$ is a $\mathbb{T}^2$-basic action, $\phi_2$ has two closed orbits of the form $\{ \zeta \} \times \text{S}^1$ and $\{ j_1' \star \zeta \} \times \text{S}^1$, for some $\zeta =\alpha e_1 + \beta e_{p+1} \neq \pm e_{p+1}$, with $ \alpha^2 + \beta^2=1 $.

\textbf{Claim 1}: There exist $z_0 \in \{e_{p+1} \} \times \text{S}^1$ and $z_0' \in \{-e_{p+1} \} \times \text{S}^1$ such that their isotropy Lie algebras with respect to $\phi_2$ are generated by the element  $X_0 = \begin{bmatrix}
     && 1 \\
     &0&\\
     1&&
 \end{bmatrix} \in \mathfrak{so}(1,2)$.

\textbf{Claim 2}: As $\theta \rightarrow \infty$, $m_\theta \star z_0$ approaches a fixed point of $\mathcal{M}(1,2)$ in a closed $\phi_2$-orbit. As $\theta \rightarrow -\infty$, $m_\theta \star z_0$ approaches a fixed point of $\mathcal{M}(1,2)$ in the other closed $\phi_2$-orbit. The same thing holds for $m_\theta \star z_0'$. 

The proof of Claims 1 and 2 is simple but a bit technical and so it is omitted; the reader can find it in \cite{Lentas_thesis}.
Let $\mathcal{S}_0$ be the $\mathcal{M}(1,2)$-orbit of $z_0$, namely
\[ \mathcal{S}_0 = \{ m(\theta) \star z_0 \}\]
By Claim 2, $\bar{\mathcal{S}_0}$ comprises $\mathcal{S}_0$ and two fixed points of $\mathcal{M}(1,2)$, one in each closed $\phi_2$-orbit. Moreover, 
it is not difficult to see that $\mathcal{S}_0$ is transverse to the $\text{SO}(2)$-orbits. Let 
\[\mathcal{S}'_0 =\{ m(\theta) \star z_0' \} \]
By Claim 2, as $\theta \rightarrow + \infty$, $m(\theta)\star z_0'$ converges to a fixed point of $\mathcal{M}(1,2)$ in one of the closed orbits. We can assume that this point is in the closure of $\mathcal{S}_0$. If it is not, we just consider $j_2\star z_0'$ instead of $z_0'$, which is easily seen to also satisfy Claims 1 and 2. Note that in the closed orbits, $j_2$ interchanges the two fixed points of $\mathcal{M}(1,2)$. By the $j_1'$ invariance of $\mathcal{S}'_0$, as $\theta \rightarrow - \infty$, $m(\theta)\star z_0'$ converges to the point in the intersection of the other closed $\phi_2$-orbit and the closure of $\mathcal{S}_0$. Finally, let
\begin{equation} \label{T2 S submnfld}
 \mathcal{S} = \overline{\mathcal{S}_0} \cup \overline{\mathcal{S}'_0} 
 \end{equation}
It can be shown that $\mathcal{S}$ is an analytic submanifold by using the local form of the infinitesimal generators of the induced $\text{SL}_2(\mathbb{R})$ action around the $\mathcal{M}(1,2)$ fixed points, see \cite[Theorem 2]{SchneiderSL2axns}. The action of $\mathcal{M}(1,2)$ on $\mathcal{S}$ induces an analytic flow $\Phi_\theta$. The Jacobian of $\Phi_\theta$ at the fixed points is non zero, see \cite[Proposition 6.1]{SchneiderSL2axns}. Since $\mathcal{S}$ is a closed, one dimensional manifold, one of the fixed points of $\Phi_\theta $ is attracting and the other is repelling. We also define a function $f_0:\mathcal{S} \rightarrow \mathbb{R}\text{P}^1$ by 
\begin{itemize}
    \item $f_0(z_0) := [0:1]$ and $ f_0(z_0'):= [0:1]$.
    \item $f_0 \left( \Phi_\theta (z) \right) := \left[ \tanh{(\theta)} : 1 \right] $, for $\theta \in \mathbb{R}$ and $z = z_0$ or $z_0'$.
    \item $f_0(z_1):=[1:1]$ and $f_0(z_2):=[1:-1]$, where $z_1$, respectively $z_2$, is the attracting, respectively repelling, fixed points of $\Phi_\theta$.
\end{itemize}
By a proof similar to that of \cite[Lemma 2.3]{Lentas_SOpq_I}, it can be shown that $f_0$ is analytic and that $\mathcal{S}$, $\Phi_\theta$, and $f_0$ satisfy the hypotheses of Lemma \ref{S f0 Phi -> (phi,f)}. Therefore, applying that lemma, we get a pair $(\phi, f)$ of an analytic $N^+$ action, $\phi$, on $\mathcal{F}$ and an analytic function $f:\mathcal{F} \rightarrow \mathbb{R}\text{P}^2$ satisfying relations (B1)-(B2). Then, by Lemma \ref{(phi , f) -> action lemma} we obtain an analytic $G$ action on $\text{S}^p \times \text{S}^1$ extending the standard $K$ action. We will refer to such actions as \textit{actions from the $(p,2)$-basic construction} (I). We note also, that since $N^+$ and $\text{SO}^\circ(1,2)^+$ are isomorphic, $\phi$ induces an action of $\text{SO}^\circ(1,2)^+$ on $\mathcal{F}$. That action, by Remark \ref{G axns <-> (phi,f)}, is $\phi_2$. 

\begin{theorem} \label{bsc constr (I) thm}
Let $p\geq3$. Two analytic $\textup{SO}^\circ(p,2)$ actions on $\textup{S}^p \times \textup{S}^1$ from the $(p,2)$-basic construction $\textup{(I)}$ are analytically isomorphic if and only if the corresponding $\mathbb{T}^2$-basic actions of $\textup{SO}^\circ(1,2)^+$ on $\mathbb{T}^2$, see Definitions \ref{SO(1,2)+ def} and \ref{T2-basic axns}, are analytically isomorphic. 
\end{theorem}
\begin{proof}
The only if part is easy. For the other implication,  let \[\mathcal{F} = \left\{ \alpha e_1 + \beta e_{p+1} : \alpha^2 + \beta^2=1 \right\} \times \text{S}^1\] 
be the fixed point set of $H$. For an analytic $\textup{SO}^\circ(p,2)$ action on $\textup{S}^p \times \textup{S}^1$ from the $(p,2)$-basic construction (I), by Remark \ref{SO(p,2) axn -> (phi,f)}, we get a pair $(\phi^1,f_{\phi^1})$ of an analytic $N^+$-action on $\mathcal{F}$ and an analytic function $ f_{\phi^1}:\mathcal{F} \rightarrow \mathbb{R}\text{P}^2$ satisfying (B1)-(B2), as well as a triplet $(\mathcal{S}^1,\Phi^1_\theta,f_0^1)$ of a closed one dimensional submanifold $\mathcal{S}^1 \subseteq \mathcal{F}$, a flow $\Phi^1_\theta$ on $\mathcal{S}^1$, and a function $f_0^1:\mathcal{S}^1 \rightarrow \mathbb{R}\text{P}^1$.  
Suppose we have a second $\textup{SO}^\circ(p,2)$ action on $\textup{S}^p \times \textup{S}^1$ from the $(p,2)$-basic construction (I), to which similarly correspond a pair $(\phi^2,f_{\phi^2})$ and a triplet $(\mathcal{S}^2,\Phi_\theta^2,f_0^2)$. Assume that the corresponding $\mathbb{T}^2$-basic actions of $\textup{SO}^\circ(1,2)^+$ on $\mathbb{T}^2$, say $\phi_1'$ and $\phi_2'$, are isomorphic. Identifying $\mathbb{T}^2$ with $\mathcal{F}$ via the function $\psi$ of (\ref{T2 and F isomor}), we get two actions of $\textup{SO}^\circ(1,2)^+$ on $\mathcal{F}$ which are isomorphic. Finally, since $\text{SO}^\circ(1,2)^+$ is isomorphic to $N^+$, we get two $N^+$-actions, $\phi_1''$ and $\phi_2''$, which are isomorphic. However, note that by Remark \ref{G axns <-> (phi,f)}, $\phi_1''=\phi_1$ and $\phi_2''=\phi_2$.

Hence, $\phi_1$ and $\phi_2$ are isomorphic. Let $\psi'$ be the isomorphism between them. Namely, $\psi' : \mathcal{F} \rightarrow \mathcal{F}$ is an analytic isomorphism such that, $\psi' \left(\phi_1''(n,p) \right) = \phi_2'' \left(n,\psi'(p) \right)$, for $n \in N^+$ and $p \in \mathcal{F}$. Then $\psi'(\mathcal{S}^1)=\mathcal{S}^2 $ and the flows $\Phi^1_\theta$ and $\Phi^2_\theta$ are isomorphic. Then, as in the proof of \cite[Theorem 2.1]{Lentas_SOpq_I}, it can be shown that 
\[ f_0^2 \circ \psi' = f_0^1\]
Since $\mathcal{S}^1$ and $\mathcal{S}^2$ are transverse to every $\text{SO}(2)$ orbit, and $\psi'$, $f_{\phi^1}$, and $f_{\phi^2}$ are $\text{SO}(2)$ equivariant, $f_{\phi^1}$, respectively $f_{\phi^2}$, is determined by $f_0^1$, respectively $f_0^2$. See also (\ref{f0 eq}). Therefore,
\[f_{\phi^2} \circ \psi = f_{\phi^1}\]
The proof then concludes similarly to the proof of \cite[Theorem 2.1]{Lentas_SOpq_I}
\end{proof}

We note again, see Remark \ref{Schneider T2}, that analytic $\text{SL}_2(\mathbb{R})$ action on $\mathbb{T}^2$ are classified in \cite{SchneiderSL2axns}, therefore, since $\text{SO}^\circ(1,2) \simeq \text{PSL}_2(\mathbb{R})$, the $\text{SO}^\circ(1,2)$ actions on $\mathbb{T}^2$ are classified. Two $\textup{SO}^\circ(1,2)^+$ are analytically isomorphic if and only if the restricted  $\text{SO}^\circ(1,2)$ are isomorphic via an isomorphism that commutes with the action of $j_1$.

\subsubsection{$(p,2)$-basic construction (II)} \label{(p,2)-basic const II}
Let $\phi_1$ be an $ \text{S}^2$-basic action of $\text{SO}^\circ(1,2)$ on $\text{S}^2$, see Definition \ref{S2 basic axn no j1}. We will show how to construct an $\text{SO}^\circ(p,2)$ action on $\text{S}^{p-1} \times \text{S}^2$ extending the standard $K$ action starting from $\phi_1$. Let 
\[ \mathcal{F} = \left\{ \pm e_1 \right\} \times \text{S}^2\]
be the fixed point set of $H$ in the standard $K$ action on $\text{S}^{p-1} \times \text{S}^2$ and $\mathcal{F}_1 = \left\{  e_1 \right\} \times \text{S}^2$, $\mathcal{F}_2 = \left\{ - e_1 \right\} \times \text{S}^2$. Let $\text{pr}_2 : \text{S}^{p-1} \times \text{S}^2 \rightarrow \text{S}^2$ be the projection to the second factor, and
\begin{equation*}
    \psi = \left(\text{pr}_2 \Big|_{\mathcal{F}_1} \right) ^{-1}
\end{equation*}
Then, $\psi: \text{S}^2 \rightarrow \mathcal{F}_1$ is an analytic isomorphism and via $\psi$, we get an analytic action, $\phi_2$, of $\text{SO}^\circ(1,2)$ on $\mathcal{F}_1$. Consider the points $(e_1, \pm \epsilon_3) \in \mathcal{F}_1$ which are fixed by $\text{SO}(2)$ and let $\mathcal{S}$ be the union of the closures of their orbits under the action of $\mathcal{M}(1,2)$. Recall that $\text{SL}_2(\mathbb{R})$ covers $\text{SO}^\circ(1,2) \simeq \text{PSL}_2(\mathbb{R})$ and hence, $\phi_1$ induces an $\text{SL}_2(\mathbb{R})$ action on $\mathcal{F}_1 \simeq \text{S}^2$. By Schneider's description of the $\text{SL}_2(\mathbb{R})$ actions in \cite{SchneiderSL2axns}, $\mathcal{S}$ is a one dimensional submanifold which is transverse to each $\text{SO}(2)$-orbit. It is also easy to see that $\mathcal{S}$ is $j_2$-invariant. The action of $\mathcal{M}(1,2) \leq \text{SO}^\circ(1,2)$ induces a flow $\Phi_\theta$ on $\mathcal{S}$. Since $\phi_1$ is an $ \text{S}^2$-basic action, $\phi_2$ has only one closed orbit on $\mathcal{F}_1$, and therefore $\Phi_\theta$ has two fixed points, see \cite[p. 521]{SchneiderSL2axns}. Define a function $f_0: \mathcal{S} \rightarrow \mathbb{R}\text{P}^1$ by
\begin{itemize}
    \item $f_0(e_1, \pm \epsilon_3): = [1:0]$
    \item $f_0 \left( \Phi_\theta (z) \right) := \left[  \, \tanh{(\theta)}  :  \, 1  \right] $, for $z_0 =  (e_1,\epsilon_3)$ or $(e_1, -\epsilon_3)$
    \item $f_0(z_1):=[1:1]$ and $f_0(z_2):=[-1:1]$, where $z_1$, respectively $z_2$, is the attracting, respectively repelling, fixed point of $\Phi_\theta$.
\end{itemize}
As in \cite[Lemma 2.3]{Lentas_SOpq_I}, it can be shown the $f_0$ is analytic and similarly to the $(p,2)$-basic construction (I), $\mathcal{S}$, $\Phi_\theta$, and $f_0$ satisfy the hypotheses of Lemma \ref{S f0 Phi -> (phi,f)}. Hence, we get an analytic pair $(\phi,f)$ of an $N^+$ action on $\mathcal{F}$ and a function $f:\mathcal{F} \rightarrow \mathbb{R}\text{P}^2$ that satisfy relations (B1)-(B2). Note that, since $\text{SO}^\circ(1,2) \simeq N$, the restriction of $\phi$ to $N$ induces an $\text{SO}^\circ(1,2)$ action on $\mathcal{F}_1$. That action is exactly $\phi_2$, see Remark \ref{G axns <-> (phi,f)}. Using the pair $(\phi,f)$ and similar arguments as in \cite[Section 2]{Lentas_SOpq_I}, we get an analytic action of $\text{SO}^\circ(p,2)$ on $\text{S}^{p-1} \times \text{S}^2$ extending the standard $K$ action.

\begin{theorem} \label{bsc constr (II) thm}
Let $p\geq3$. Two analytic $\textup{SO}^\circ(p,2)$ actions on $\textup{S}^{p-1} \times \textup{S}^2$ from the $(p,2)$-basic construction $\textup{(II)}$ are analytically isomorphic if and only if the corresponding $ \textup{S}^2$-basic actions of $\textup{SO}^\circ(1,2)$ actions on $\textup{S}^2$, see Definition \ref{S2 basic axn no j1}, are isomorphic.
\end{theorem}
The proof is similar to that of Theorem \ref{bsc constr (I) thm}. We note again, that the $\textup{SO}^\circ(1,2)$ actions on $\text{S}^2$ are classified by the results in \cite{SchneiderSL2axns}, see Remark \ref{Schneider S2}.

\subsubsection{$(p,2)$-basic construction (III)} \label{(p,2)-basic const III}
 Let $\phi_1$ be an $(\text{S}^2)^+$-basic action of $\text{SO}^\circ(1,2)^+$ on $\text{S}^2$, see Definition \ref{basic S2 axn + j1}. Starting from $\phi_1$, we will construct an $\text{SO}^\circ(p,2)$ action on $\text{S}^{p+1}$ extending the standard $K$ action. These are the kind of actions Uchida studied in \cite{Uchida_SOpqII} in the smooth setting. In the standard $K$ action, the fixed point set of $H$ is 
\[ \mathcal{F} = \{ \alpha \, e_1 + \beta \, e_{p+1} + \gamma \, e_{p+1} \, : \, \alpha^2 + \beta^2 + \gamma^2 =1 \}\]
Let $\psi : \text{S}^2 \rightarrow \mathcal{F}$ be defined by
\[ (\alpha \, e_1 + \beta \, e_{p+1} + \gamma \, e_{3}) \mapsto \alpha \, e_1 + \beta \, e_{p+1} + \gamma \, e_{p+1}\]
Then, $\psi$ is an analytic isomorphism. Via $\psi$, we get an action, $\phi_2$, of $\text{SO}^\circ(1,2)^+$ on $\mathcal{F}$. Let $\mathcal{S}_1$ be the $\mathcal{M}(1,2)$-orbit of the point $e_1 \in \mathcal{F}$. It is easy to see that $\mathcal{S}_1$ is $j_2$-invariant since $e_1$ is fixed by $j_2$. Now, $\text{SL}_2(\mathbb{R})$ covers $\text{SO}^\circ(1,2) \simeq \text{PSL}_2(\mathbb{R})$ and hence, we get an $\text{SL}_2(\mathbb{R})$ action on $\mathcal{F}_1 \simeq \text{S}^2$ induced by $\phi_2$. By Schneider's description of the $\text{SL}_2(\mathbb{R})$ action on $\text{S}^2$, see \cite{SchneiderSL2axns}, we can find a point $z=\alpha \, e_{p+1}+ \beta \, e_{p+2} \in \text{S}^2$ such that
\[ \lim_{\theta \rightarrow +\infty} m(\theta) \star z =\lim_{\theta \rightarrow +\infty} m(\theta) \star e_1 \]
for $m(\theta) \in \mathcal{M}(1,2)$. Let $\mathcal{S}_2$ be the $\mathcal{M}(1,2)$-orbit of $z$. Finally, consider \[\mathcal{S}= \mathcal{S}_1 \cup j_1'(\mathcal{S}_1) \cup \mathcal{S}_2 \cup j_2 \left(\mathcal{S}_2 \right)\]
Then, $\mathcal{S}$ is a one dimensional submanifold of $\text{S}^2$, it is $j_1'$ and $j_2$-invariant, and it is transverse to each $\text{SO}(2)$-orbit, see \cite{SchneiderSL2axns}. The action of $\mathcal{M}(1,2)$ on $\mathcal{S}$ induces an analytic flow, $\Phi_\theta$, on $\mathcal{S}$ . Since $\phi_1$ is an $(\text{S}^2)^+$-basic action, $\phi_2$ has two closed orbits in $\mathcal{F}$ and hence, $\Phi_\theta$ has four fixed points on $\mathcal{S}$, see also \cite{SchneiderSL2axns}. Then, we can define a function $f_0:\mathcal{S}\rightarrow \mathbb{R}\text{P}^1$ by
\begin{itemize}
    \item $f_0(\pm e_1):= [1:0]$, $f_0(z):= [0:1]$, and $f(j_2 \star z) := [0:1]$
    \item $f_0 \left( \Phi_\theta (e_1) \right) := \left[ \text{cosh}(\theta)  : \text{sinh}(\theta) \right]$ ,  $f_0 \left( \Phi_\theta (z) \right) := \left[ \text{sinh}(\theta)  : \text{cosh}(\theta) \right]$ and in the rest of $\mathcal{S}$ $f_0$ is defined so that it satisfies $j_1' \star \Phi_\theta (z) = \Phi_{- \theta} (j_1' \star z)$ and $j_2 \star \Phi_\theta (z) = \Phi_{- \theta} (j_2 \star z)$ 
    \item $f_0(\zeta):= [1:1]$, respectively $[1:-1]$, if $\zeta \in \mathcal{S}$ is an attracting, respectively repelling, fixed point of $\Phi_\theta$.
\end{itemize}

 Similarly to \cite[Lemma 2.3]{Lentas_SOpq_I}, it can be shown the $f_0$ is analytic. Consequently, analogously to the $(p,2)$-basic construction (I), $\mathcal{S}$, $\Phi_\theta$, and $f_0$ satisfy the hypotheses of Lemma \ref{S f0 Phi -> (phi,f)}, which gives an analytic pair $(\phi,f)$ of an $N^+$ action on $\mathcal{F}$ and a function $f:\mathcal{F} \rightarrow \mathbb{R}\text{P}^2$ that satisfy relations (B1)-(B2). Note that, since $\text{SO}^\circ(1,2)$ and $ N$ are isomorphic, $\phi$ induces an $\text{SO}^\circ(1,2)$ action on $\mathcal{F}$. That action is $\phi_2$, see also Remark \ref{G axns <-> (phi,f)}. Using the pair $(\phi,f)$ and analogous arguments as in \cite[Section 2]{Lentas_SOpq_I}, we get an analytic action of $\text{SO}^\circ(p,2)$ on $\text{S}^{p-1} \times \text{S}^2$ extending the standard $K$ action. The following theorem strengthens \cite[Theorem 4.12]{Uchida_SOpqII} in the analytic setting.

\begin{theorem} \label{bsc constr (III) thm} \textup{(See also \cite[Theorem 4.12]{Uchida_SOpqII})} 
Let $p\geq3$. Analytic $\textup{SO}^\circ(p,2)$ actions on $\textup{S}^{p+1}$ from the $(p,2)$-basic construction $\textup{(III)}$ are analytically isomorphic if and only if the corresponding $\left( \textup{S}^2 \right)^+$-basic actions of $\textup{SO}^\circ(1,2)^+$ actions on $\textup{S}^2$, see \textup{(\ref{SO(1,2)+ def})} and Definition \ref{basic S2 axn + j1}, are isomorphic.
\end{theorem}
 The proof is similar to that of Theorem \ref{bsc constr (I) thm}. We note, that the $\textup{SO}^\circ(1,2)^+$ actions on $\text{S}^2$ are classified by the results in \cite{SchneiderSL2axns}, see Remark \ref{schneider S2+}.

\subsection{Classification}

\begin{theorem}
Suppose $\textup{SO}^\circ(p,2)$, $p \geq 3$, acts analytically on a  closed, connected manifold $M$ of dimension $p+1$. Consider $\textup{SO}(p) \simeq \textup{SO}(p) \times \{1\} \leq \textup{SO}(p) \times \textup{SO}(2) \leq \textup{SO}^\circ(p,2)$ and $\textup{SO}(2)$ similarly. For $p>3$, we have:
\begin{itemize}
    \item If only $\textup{SO}(p)$ has a fixed point, then $M$ is equivariantly covered by $\textup{S}^p \times \textup{S}^{1}$, where the action of $\textup{SO}^\circ(p,2)$ on $\textup{S}^p \times \textup{S}^{1}$ is one from the $(p,2)$-basic construction \textup{(I)}, see Section \ref{(p,2)-basic const I}.
    \item If only $\textup{SO}(2)$ has a fixed point, then $M$ is equivariantly covered by $\textup{S}^{p-1} \times \textup{S}^{2}$, where the action of $\textup{SO}^\circ(p,2)$ on $\textup{S}^{p-1} \times \textup{S}^2$ is one from the $(p,2)$-basic construction \textup{(II)}, see Section \ref{(p,2)-basic const II}.
    \item If both $\textup{SO}(p)$ and $\textup{SO}(2)$ have a fixed point, then $M$ is equivariantly covered by $\textup{S}^{p+1}$, where the action of $\textup{SO}^\circ(p,2)$ on $\textup{S}^{p+1}$ is one from the $(p,2)$-basic construction \textup{(III)}, see Section \ref{(p,2)-basic const III}.
    \item If neither $\textup{SO}(p)$ nor $\textup{SO}(2)$ have a fixed point, then $M$ is equivariantly covered by $\textup{SO}^\circ(p,2) \times_P \textup{S}^1$ with the left $\textup{SO}^\circ(p,2)$ action, where $P \leq \textup{SO}^\circ(p,2)$ is a maximal parabolic subgroup isomorphic to the stabiliser of an isotropic line in the standard representation of $\textup{SO}^\circ(p,2)$ on $\mathbb{R}^{p+2}$. If $P=M_P A_P N_P$ is a Langlands decomposition of $P$, $P$ acts on $\textup{S}^1$ by a flow via $A_P$, see \cite[(21)]{Lentas_SOpq_I}.
\end{itemize}
For $\textup{SO}^\circ(3,2)$, the above applies if the action is not transitive and covered by $\displaystyle \textup{SO}^\circ(3,2) \big/ P_{\textup{min}}$, where $P_{\textup{min}} \leq \textup{SO}^\circ(3,2)$ is a minimal parabolic subgroup.
\end{theorem}
\begin{proof}
Suppose $G=\textup{SO}^\circ(p,2)$  acts analytically on a manifold $M$ of dimension $p+1$. In the case of $\text{SO}^\circ(3,2)$ the homogeneous space $\textup{SO}^\circ(3,2) \big/ P_{\textup{min}}$, where $P_{\textup{min}}$ is a minimal parabolic subgroup of $\text{SO}^\circ(3,2)$ has dimension $4 = p+1$. So, suppose that either $p>3$ or in the case of $p=3$ that the action is not covered by $\textup{SO}^\circ(3,2) \big/ P_{\textup{min}}$. Then, the methods of \cite[Section 3]{Lentas_SOpq_I} can be used to show that the fixed point set of $H$, $\mathcal{F}$, is non-empty and 2-dimensional. We will split the proof in three cases depending on the existence of fixed points of $\text{SO}(p)$ or $\text{SO}(2)$.

Similarly to Remark \ref{SO(p,2) axn -> (phi,f)}, from this action of $G$ we get the restricted action, $\phi_1$, of $N^+$ on $\mathcal{F}$ and a function $f_1:\mathcal{F} \rightarrow \mathbb{R}\text{P}^2$, by Lemma \ref{SOp,2 lie alg lemma}. Additionally, $\phi_1$ and $f_1$ satisfy properties (B1)-(B2) from Remark \ref{SO(p,2) axn -> (phi,f)}, see \cite{Uchida_SOpqII}.  Let $\mathcal{R} = \{ [a:b:c] \in \mathbb{R}\text{P}^2 : c=0 \} \subseteq \mathbb{R}\text{P}^2 $ and consider $S = f_1^{-1}(\mathcal{R}) \subset \mathcal{F}$.

\textbf{Case I:} Suppose that $\text{SO}(p)$ has a fixed point and that $\text{SO}(2)$ does not have any fixed points. Since $N \leq G$ acts on $\mathcal{F}$,  and $ N \simeq \text{PSL}_2(\mathbb{R})$, see the beginning of this section, $\text{SL}_2(\mathbb{R})$ also acts on $\mathcal{F}$. However, $\text{SO}(2)$ does not have fixed points, hence by \cite{SchneiderSL2axns}, $\mathcal{F}$ has to be diffeomorphic to a torus. As in \cite[Section 4.2]{Lentas_SOpq_I}, we can find an analytic isomorphism 
\[ \psi : S \rightarrow \text{S}^1 \quad \text{s.t.} \quad \psi \circ j_1 = j_1 \circ \psi \]
where $j_1$ acts on $\text{S}^1$ by reflection with respect to the $y$-axis. Suppose the image of $\psi$ is the first factor of $\mathbb{T}^2 = \text{S}^1 \times \text{S}^1$.
Then, define:
\[ \Psi : \mathcal{F} \rightarrow \mathbb{T}^2 \]
by
\[ R_\theta z \mapsto R_\theta \psi (z) \]
for $z \in S$ and $R_\theta \in \text{SO}(2)$, where $\text{SO}(2)$ acts on the second factor of $\mathbb{T}^2$ in the standard way. Recall that $S$ is transverse to all the $\text{SO}(2)$ orbits. Now, if we identify $\mathbb{T}^2$ with the set
\[\hat{\mathcal{F}} = \left\{ \alpha e_1 + \beta e_{p+1} : \alpha^2 + \beta^2=1 \right\} \times \text{S}^1\]
in the obvious way, then, via $\Psi$ and $(\phi_1,f_1)$, we get a pair $(\phi , f)$ of an analytic $N^+$ action on $\hat{\mathcal{F}}$ and an analytic function $f:\hat{\mathcal{F}} \rightarrow \mathbb{R}\text{P}^2$ satisfying (B1)-(B2) of Lemma \ref{(phi , f) -> action lemma}and hence an analytic action of $G$ on $\text{S}^p \times \text{S}^1$ that extends the standard $K$ action. Note that this is an action from the basic construction corresponding to $\phi'$, where $\phi'$ is $\mathbb{T}^2$-basic action of $\text{SO}^\circ(1,2)^+$ on $\mathbb{T}^2$ induced by $\Psi$ and the isomorphism between $N^+$ and $\text{SO}^\circ(1,2)^+$, see (\ref{N^+ SO(1,2)+ isom}). For $p \in \text{S}^p \times \text{S}^1$, there exist $k \in K$ and $Y \in \hat{\mathcal{F}}$ such that $k \star Y = p$. Then, we define
\begin{align*}
    F: & \text{S}^p \times \text{S}^1  \rightarrow M \\
    & k\star Y  \mapsto k \star \Psi(Y)
\end{align*}
As in \cite[Section 4.2]{Lentas_SOpq_I} it can be shown that all the resulting maps are analytic and $F$ is a covering map.

\textbf{Case II:} In the case that neither $\text{SO}(p)$ nor $\text{SO}(q)$ has a fixed point, looking at $S \subset \mathcal{F}$ we get that
\[ f_1 \Big|_J = [1: \pm 1:0]\]
where $J \subseteq S$ is a nontrivial open subset. Without loss of generality, we may assume that the value of $f_1$ on $J$ is $[1:1:0]$. Since, $f$ is analytic,
\[ f_1 \Big|_J \equiv [1:1:0]\]
That means that the isotropy group of every point on $S$ is contained in or equal to $P$, where $P \leq G$ is a maximal parabolic. As in \cite[Section 4.1]{Lentas_SOpq_I}, it can be shown that such an action is covered by an associated bundle $G \times_P \text{S}^1$. Note that the absence of  $\text{SO}(p)$ and $\text{SO}(2)$ fixed points is equivalent to the existence of a nullcone orbit.

\textbf{Case III:} Now, suppose that $\text{SO}(2)$ has a fixed point. Recall that via $N \simeq \text{PSL}_2(\mathbb{R})$, $\text{SL}_2(\mathbb{R})$ also acts on $\mathcal{F}$. According to Schneider's results, the connected components of $\mathcal{F}$ are diffeomorphic to either $\mathbb{R}\text{P}^2$ or $\text{S}^2$, see \cite{SchneiderSL2axns}. Suppose a component is diffeomorphic to $\mathbb{R}\text{P}^2$
and denote this component by $\mathcal{F}_1$. Assume that $j_1$ does not preserve this component. By Schneider's description of the $\text{SL}_2(\mathbb{R})$ actions on $\mathbb{R}\text{P}^2$, there exists an action of $\text{SL}_2(\mathbb{R})$ and hence an action, $\phi_2$, of $N$ on $\text{S}^2$ that covers the action $\phi_1$ of $N$ on $\mathcal{F}_1$ equivariantly. Let $\pi : \text{S}^2 \rightarrow \mathcal{F}_1$ be that $N$-equivariant covering map. Let $\rho_1$ be the obvious inclusion
\[\rho_1 :\text{S}^2 \rightarrow \left\{ e_{1} \right\} \times \text{S}^2\]
Via $\rho_1$, $\phi_2$ induces an action $\phi'$ of $N$ on $\left\{ e_{1} \right\} \times \text{S}^2 \subseteq \text{S}^{p-1} \times \text{S}^2$. Now, on $\left\{ -e_{1} \right\} \times \text{S}^2$ consider the action of $N$ defined by
\[ n \star (j_1v) = j_1 \left( j_2nj_2 \star v \right)\]
for $n \in N$ and $v \in  \left\{ e_{1} \right\} \times \text{S}^2$. Then, $ \left\{- e_{1} \right\} \times \text{S}^2$ covers $N$-equivariantly the component $j_1 \left( \mathcal{F}_1 \right)$ and we get an action $\phi$ of $N^+$ on $ \left\{ \pm e_{1} \right\} \times \text{S}^2$, extending $\phi'$. Now, let
\[ f_2 :  \left\{ \pm e_{1} \right\} \times \text{S}^2 \rightarrow \mathbb{R}\text{P}^2\]
be defined by $f_2 = f_1 \circ \pi \circ \rho_1^{-1} $ on $\left\{ e_{1} \right\} \times \text{S}^2$ and by $f_2 = f_1 \circ \pi \circ \rho_1^{-1}\circ j_1 $ on $ \left\{- e_{1} \right\} \times \text{S}^2$. Then, 
$\phi$ and $f_2$ satisfy (B1)-(B2) of Lemma \ref{(phi , f) -> action lemma}. As in \cite[Section 4.2]{Lentas_SOpq_I} it can be shown that there exists a $G$ action on $\text{S}^{p-1} \times \text{S}^2$ from the $(p,2)$-basic construction (II), which covers the action on $M$, $G$-equivariantly. A similar argument applies when the component of $\mathcal{F}$ is diffeomorphic to $\text{S}^2$ and $j_1$ does not preserve this component.

 Now, suppose again that a component of $\mathcal{F}$ is diffeomorphic to $\mathbb{R}\text{P}^2$ and denote this component by $\mathcal{F}_1$, but now we assume that $j_1$ preserves $\mathcal{F}_1$. Now, $S$ is a 1-dimensional closed submanifold of $\mathcal{F}_1$. Let $S_0$ be a connected component of $S$ which contains a fixed point of $\text{SO}(2)$, say $Y$. We claim that $j_1$ preserves $S_0$. Indeed, otherwise $j_1Y$ would also be an $\text{SO}(2)$ fixed point and there would be at least two points fixed by $\text{SO}(2)$. But that is impossible for an $\text{SL}_2(\mathbb{R})$ action on $\mathbb{R}\text{P}^2$, see \cite{SchneiderSL2axns}. Therefore $j_1$ acts on $S_0$. Note that, being an involution on a manifold diffeomorphic to $\text{S}^1$, $j_1$ has either none or two fixed points. Let $Y \in \mathcal{F}_1$ be the fixed point of $\text{SO}(2)$. Since $j_1$ commutes with $\text{SO}(2)$ in $G$,  $j_1Y$ is also fixed by $\text{SO}(2)$ and hence it must equal $Y$. So, $j_1$ has at least one fixed point on $S_0$ and therefore it has exactly two. Suppose $j_1$ acts on $\mathbb{R}\text{P}^1$ by 
\[ j_1 \star [a:b] = [-a:b]\]

A slight variation in the proof of \cite[Lemma 4.4]{Lentas_SOpq_I} shows that there exists an analytic isomorphism 
\[ \psi : \mathbb{R}\text{P}^1 \rightarrow S_0\]
which is $j_1$ equivariant. Define an action of $\mathcal{M}(p,2)$ on $\mathbb{R}\text{P}^1$ by
\[ m(\theta) \star \psi^{-1}(x) : = \psi^{-1} (m(\theta) \star x) \]
and a function $f_2 : \mathbb{R}\text{P}^1 \rightarrow \mathbb{R}\text{P}^1$ by
\[ f_2 := f_1 \circ \psi\]
Imposing $\text{SO}(2)$-equivariance, we can extend $\psi$ to a map $\hat{\psi}$,
\[ \hat{\psi} : \mathbb{R}\text{P}^2 \rightarrow \mathcal{F}_1\]
and using $f_2$ and the action of $\mathcal{M}(p,2)$ we can define an action $\hat{\phi}$ of $N^+$ on $\mathbb{R}\text{P}^2$. The proofs are completely analogous to the arguments in \cite[Section 2]{Lentas_SOpq_I}, using (\ref{N decomp}). Then, $\hat{\psi}$ is $N^+$-equivariant.

Now, as before we can find an action, $\phi'$, of $N$ on $\text{S}^2$ such that it covers the $N$ action $\hat{\phi}$ on $\mathbb{R}\text{P}^2$. On $\text{S}^2$ we consider the action of $j_1$ by
\[ j_1 \star (a \, e_1 +b \, e_2+c \, e_3) = -a \, e_1 +b \, e_2+c \, e_3\]
Then the projection to $\mathbb{R}\text{P}^2$ is also $j_1$ equivariant. Composing with $\hat{\psi}$, we get a covering map
\[ \tilde{\psi} : \text{S}^2 \rightarrow \mathcal{F}_1\]
which is $N^+$-equivariant. Composing this map with $f_1$, we get a map
\[ f_2'=f_1 \circ \tilde{\psi} : \text{S}^2 \rightarrow \mathbb{R}\text{P}^2\]
Let
\[ \hat{\mathcal{F}} = \{ \alpha \, e_1 + \beta \, e_{p+1} + \gamma \, e_{p+1} \, : \, \alpha^2 + \beta^2 + \gamma^2 =1 \} \subseteq \text{S}^{p+1} \]
and let $\rho_2 : \text{S}^2 \rightarrow \hat{\mathcal{F}}$ be the isomorphism defined by
\[ \rho_2 \left( \alpha \, e_1 + \beta \, e_{p+1} + \gamma \, e_{3} \right) = \alpha \, e_1 + \beta \, e_{p+1} + \gamma \, e_{p+1} \]
Then, via $\rho_2$ we get an action, $\phi$, of $N^+$ on $\hat{\mathcal{F}}$. Finally, consider the function
\[f = f'_2 \circ \rho_2^{-1}\]
As in the first part of Case III, $(\phi,f)$ satisfy relations (B1)-(B2) from Lemma \ref{(phi , f) -> action lemma}, and therefore, there exists a $G$ action on $\text{S}^{p+1}$ from the $(p,2)$-basic construction (III), such that $\text{S}^{p+1}$ equivariantly covers $M$. Similar arguments apply when the component of $\mathcal{F}$ is diffeomorphic to $\text{S}^2$ and $j_1$ preserves this component.
\end{proof}
\printbibliography
\end{document}